\documentclass{article}

\usepackage{amsmath, amsthm, amsfonts, amssymb, hyperref, enumerate, esint, fullpage,dsfont}
\usepackage{mathtools}

\usepackage[dvipsnames]{xcolor}

\usepackage{todonotes}

\newtheorem{theorem}{Theorem}
\newtheorem{lemma}{Lemma}
\newtheorem{prop}{Proposition}

\theoremstyle{definition}
\newtheorem{defn}{Definition}

\theoremstyle{remark}
\newtheorem{remark}{Remark}

\newcommand{\ud}{\,\mathrm{d}}
\newcommand{\R}{\mathbb{R}}
\newcommand{\Rd}{\mathbb{R}^d}
\newcommand{\Z}{\mathbb{Z}}
\newcommand{\Zd}{\mathbb{Z}^d}

\newcommand{\Td}{\mathbb{T}^d}

\newcommand{\ind}{\mathrm{1}}

\newcommand{\Sone}{\mathcal{S}^1}
\newcommand{\Stwo}{\mathcal{S}^2}

\newcommand{\norm}[1]{\left\lVert#1\right\rVert}

\begin{document}
	\title{On Adaptive Confidence Sets for the Wasserstein Distances}
	
	\author{
	Neil Deo\thanks{University of Cambridge (\texttt{and30@cam.ac.uk})}
	\and
	Thibault Randrianarisoa\thanks{Sorbonne Université (\texttt{thibault.randrianarisoa@sorbonne-université.fr})}
	}
	\date{}
	\maketitle
	\begin{center}
    	\vspace{0.9cm}
    	\textbf{Abstract}
	\end{center}
	
	%Building on recent work from J. Niles-Weed and Q. Berthet, 
	In the density estimation model, we investigate the problem of  constructing \textit{adaptive honest confidence sets} with radius measured in Wasserstein distance $W_p$, $p\geq1$, and for densities with unknown regularity measured on a Besov scale. As sampling domains, we focus on the $d-$dimensional torus $\mathbb{T}^d$, in which case $1\leq p\leq 2$, and $\mathbb{R}^d$, for which $p=1$. We identify necessary and sufficient conditions for the existence of adaptive confidence sets with diameters of the order of the regularity-dependent $W_p$-minimax estimation rate. Interestingly, it appears that the possibility of such adaptation of the diameter depends on the dimension of the underlying space. In low dimensions, $d\leq 4$, adaptation to any regularity is possible. In higher dimensions, adaptation is possible if and only if the underlying regularities belong to some interval of width at least $d/(d-4)$. This contrasts with the usual $L_p-$theory where, independently of the dimension, adaptation requires regularities to lie in a small fixed-width window. For configurations allowing these adaptive sets to exist, we explicitly construct confidence regions via the method of risk estimation, centred at adaptive estimators. Those are the first results in a statistical approach to adaptive uncertainty quantification with Wasserstein distances. Our analysis and methods extend more globally to weak losses such as Sobolev norm distances with negative smoothness indices.
	
	\tableofcontents

	\section{Introduction}

The construction of confidence sets is one of the fundamental problems of statistical inference, along with parameter estimation and hypothesis testing. Consider a model $\{P_f: f\in\mathcal{F}\}$, indexed by a family of functions $\mathcal{F}$, and observe (some quantity $n$ of) data from the true distribution $P_{f_0}$, where $f_0\in\mathcal{F}$. For most applications, having a single point estimate $\hat{f}_n$ of the true parameter $f_0$ is not enough, and one desires to evaluate its performance in terms of a loss function, that is, to know how far it lies from $f_0$. Producing a random set $C_n\subset\mathcal{F}$ from the data containing $f_0$ with a prescribed high probability $1-\alpha$ achieves this aim. In this work, we investigate the existence of \emph{adaptive honest confidence sets}. Since $f_0$ is unknown, we must insist that $C_n$ possesses the previous property not just for $f_0$, but for all $f\in\mathcal{F}$: we say that the confidence set $C_n$ is \emph{honest} if, at least for all sufficiently large $n$,
$$ \inf_{f\in\mathcal{F}}P_f(f\in C_n) \geq 1-\alpha.$$
Furthermore, we desire the diameter of the set $C_n$ to shrink in $n$ as quickly as possible; however, typically the precise speed of this shrinkage depends on aspects of the unknown density $f_0$ such as its regularity, and so we find ourselves in an adaptation problem.

We work in a \emph{density estimation model}: consider observations $X_1,\dots,X_n$ independent and identically distributed (i.i.d.) from a probability measure $P_{f_0}$ with probability density $f_0$. The sample space of the $X_i$'s will either be the $d-$dimensional torus $\mathbb{T}^d$ or $\mathbb{R}^d$. We then study procedures in a representative `two-class adaptation problem', where $f_0$ belongs to one of two classes $\mathcal{F}(r)$ and $\mathcal{F}(s)$ (to be precisely defined below), indexed by regularity parameters $r<s$, such that $\mathcal{F}(s)\subset\mathcal{F}(r)$. An adaptive honest confidence set $C_n$ should satisfy the above honest coverage condition, and also have a diameter that shrinks at the minimax estimation rate of whichever class $f_0$ belongs to (typically the rate is faster for the smaller class $\mathcal{F}(s)$). The construction of such a confidence set involves assessing the accuracy with which one can estimate $f_0$, which turns out to be more challenging than point estimation, as qualitative aspects of the parameter need to be identified. This problem has primarily been studied for $L_p$ or related distances \cite{lowNonparametricConfidenceIntervals1997, juditskyNonparametricConfidenceSet2003, caiAdaptiveConfidenceBalls2006,robinsAdaptiveNonparametricConfidence2006, hoffmannAdaptiveInferenceConfidence2011, bullAdaptiveConfidenceSets2013, carpentierHonestAdaptiveConfidence2013}. In $L_2$ loss, adaptive honest confidence sets exist only if the regularity parameters of interest lie in some `small' interval. More troublesome is the case of pointwise or $L_\infty$ loss, where no such procedures exist. This starkly contrasts the situation of adaptive estimation, where (perhaps at the cost of a logarithmic factor) it is possible to construct estimators which adapt to any regularity parameter (\cite{lepskiiProblemAdaptiveEstimation1991, donohoDensityEstimationWavelet1996}). Informally, these negative results come from the fact that, in $L_2$ loss, a related testing problem is easier (admits a faster convergence rate) than estimation, whereas for $L_\infty$ loss, the testing and estimation problems are equally difficult (\cite{hoffmannAdaptiveInferenceConfidence2011, bullAdaptiveConfidenceSets2013}). This distinction highlights how the existence of adaptive honest confidence sets depends on the geometry induced by the loss function (see \cite[Chapter 8]{gineMathematicalFoundationsInfiniteDimensional2015} for an overview of these results).

Arising from the ideas of Optimal Transport \cite{monge1781memoire, Kant42}, Wasserstein distances $W_p,\ p\geq1$, between probability measures have recently been studied in a wide array of fields such as optimization, machine learning, and statistics. For $p\geq1$, the $p-$Wasserstein distance between $\mu$ and $\nu$, probability measures on a metric space $\left(\mathcal{X}, d\right)$, is defined as
\[ W_p(\nu,\mu) \coloneqq \underset{\pi\in\Pi\left(\nu,\mu\right)}{\inf}\left(\int_{\mathcal{X}\times\mathcal{X}} d(x,y)^p d\pi(x,y)\right)^{1/p},\]
with the infimum ranging over the set $\Pi\left(\nu,\mu\right)$ of measures on $\mathcal{X}\times\mathcal{X}$ with given marginals $\nu$ and $\mu$. It quantifies the minimal cost, as measured by the metric $d$, to morph the distribution $\mu$ into $\nu$. For measures $P_f$ and $P_g$ dominated by a common measure and with densities $f$ and $g$, this also entails a distance between those densities, with $W_p(f,g)\coloneqq W_p\left(P_f,P_g\right)$.

Not only do these distances possess desirable theoretical properties (\cite{villaniOptimalTransportOld2009}), as they take into account the geometry of the underlying sample space, but recent numerical developments (\cite{COTFNT}) have led to increased use in practical applications. They therefore now play a prominent role in statistics (see, for example, the review \cite{PanaretosZemel}). The convergence of the empirical distribution in $W_p$-distance is a well-studied problem (it stretches back to \cite{Dud69}, with definitive results on limit theorems for the $\mathbb{R}$ sample space in \cite{MR1698999}; for state-of-the-art results, see \cite{fournierRateConvergenceWasserstein2015, weedSharpAsymptoticFinitesample2019}). In dimensions $d\geq3$, the convergence rate of the empirical distribution (without further structural assumptions) is $n^{-1/d}$, demonstrating that convergence in $W_p$ suffers from the curse of dimensionality.
% Beside this estimation problem, the statistical application of the Wasserstein distance has also permit the development of a theory of goodness-of-fit tests \cite{delBarrio1999, dewet2002}. 
% Further, \cite{2017arXiv171010457D} derived the minimax separation rate of tests with simple null hypotheses.
When measures have densities, as is the case in density estimation, \cite{weedEstimationSmoothDensities2019} prove that, for certain classes of densities, $W_p$ compares with Besov norms of smoothness $-1$, a classical result for the $W_1$ distance due to the Kantorovich-Rubinstein duality formula. The convergence rates they obtain for regular densities using this comparison result, which lie closer to the parametric rate $n^{-1/2}$, highlight the importance of regularity of the signal in high-dimensional settings: to some extent, the curse of dimensionality can be mitigated by smoothness.

In addition, these rates are faster than the standard $s$-smooth nonparametric convergence rate $n^{-\frac{s}{2s+d}}$ for $L_p$ loss, $1\leq p<\infty$, reflecting the fact that Wasserstein distances are weaker than $L_p$ distances. In this paper, we obtain similar quantitative improvements for testing separation rates of nonparametric statistical hypotheses. From this, on the bounded sample space $\Td$ we deduce new qualitative phenomena regarding the existence and non-existence of adaptive honest confidence sets when using the loss functions $W_p$, $1\leq p\leq 2$. Surprisingly, in dimensions $d\leq 4$ we construct confidence sets that can adapt to \emph{any} set of regularities. This contrasts significantly with the fundamental limitations of adaptive confidence sets in $L_p$. In higher dimensions $d>4$, adaptation is still possible for regularities belonging to a certain interval, which is wider than in the $L_p$ case. The reason for this phenomenon is that while both the testing and estimation rates are faster than for $L_p$, the testing rate accelerates more, leaving `more space' for adaptation to occur than in the analogous problem for $L_p$ loss. As for densities on an unbounded sample space such as $\mathbb{R}^d$, the same phenomenon occurs, though we currently only have results for the $W_1$ distance.

The paper is organized as follows. Section \ref{Section: Main Results} formalizes our problem on the potential existence of adaptive honest confidence sets, and states our main results. The construction of such sets, whenever possible, and non-existence results are presented in Section \ref{Section: Existence Results} for the bounded sample space $\Td$ and Section \ref{Section: Rd} for the unbounded sample space $\Rd$. Proofs are deferred to Appendices \ref{Section: additional torus proofs} and \ref{Section: Rd proofs}.

	\section{Main Results}\label{Section: Main Results}

\subsection{Setting and Definitions}

Initially, we assume that $f_0$ is a density on the $d$-dimensional torus, $\Td$, which may be identified with $(0,1]^d$. Our results also apply to the case of the unit cube $[0,1]^d$ (and hence any bounded rectangular subset of $\Rd$), which is the focus of \cite{weedEstimationSmoothDensities2019}; see Section \ref{Subsection: unit cube} below. For our loss function, we take the distance $W_2$; as described in Remark \ref{Remark: choice of W_2}, this distance dominates $W_p$ for $1\leq p<2$, in particular the important case of $W_1$.
Later, we consider the situation where $f_0$ is a density on the whole of $\Rd$; while a study for $W_p,p>1$ is beyond the scope of the present work, we obtain some definitive results for the loss function $W_1$ in Section \ref{Section: Rd}.

\subsubsection{Parameter Spaces}

Here we define the classes of probability densities on $\Td$ we consider; definitions for $\Rd$ are similar but deferred to Section \ref{Section: Rd}. Let $ \left\{ \phi \equiv 1, \psi_{lk} : l\geq0, 0\leq k<2^{ld}\right\}$ be an $S$-regular periodised Daubechies wavelet basis of $L_2(\Td)$; see Appendix \ref{Section: Wavelet Appendix} for further details. We denote by $\langle f,g \rangle = \int_{\Td} fg$ the usual inner product on $L_2$. For any $f\in L_p(\Td), 1\leq p<\infty$, the wavelet expansion
\begin{equation}\label{Eq: wavelet expansion}
	f = \langle f, 1 \rangle + \sum_{l\geq0}\sum_{k=0}^{2^{ld}-1} \langle f,\psi_{lk}\rangle \psi_{lk}
\end{equation}
converges in $L_p$, and if $f$ is continuous then the expansion converges uniformly on $\Td$. We write $K_j(f)$ for the projection of $f$ onto the first $j$ resolution levels, i.e.
\begin{equation}\label{Eq: j level projection definition}
	K_j(f) = \langle f,1 \rangle + \sum_{l<j}\sum_{k=0}^{2^{ld}-1}\langle f,\psi_{lk} \rangle \psi_{lk}.
\end{equation}

To define the parameter classes, we use the scale of \emph{Besov spaces,} $B^s_{pq},1\leq p,q\leq\infty, s\geq0$ as defined in Appendix \ref{Section: Wavelet Appendix}. The index $s$ should be interpreted as a smoothness or regularity parameter. Using the definition of the Besov norm (\ref{Eq: Besov norm definition}) and the embedding $\ell_q\subset\ell_{\infty}$, for $f\in B^{s}_{pq}(\Td)$ we have that
\begin{equation}\label{Eq: Besov wavelet coefficient bound}
	\left\| \langle f,\psi_{l\cdot} \rangle\right\|_p \leq \|f\|_{B^s_{pq}}2^{-l\left(s + \frac{d}{2} - \frac{d}{p}\right)}.
\end{equation}
Thus $f\in B^s_{pq}$ if its wavelet coefficients decay sufficiently fast as $l$ grows, as measured by $s$. 

The use of subsets of Besov spaces as parameter spaces in nonparametric statistics is well-established, and the scale contains several of the regularity classes usually considered in such settings: for example, the Sobolev spaces ($H^s = B^s_{22}$) and the H\"{o}lder spaces (for $s\not\in\mathbb{N}, C^s=B^s_{\infty\infty}$, and for $s\in\mathbb{N},C^s\subsetneq B^s_{\infty\infty}$). See \cite[Section 4.3]{gineMathematicalFoundationsInfiniteDimensional2015} for further discussion on this subject.

In standard loss functions such as $L_2$ or $L_{\infty}$, it is typically assumed that $f$ lies in some norm-ball in $B^s_{pq}$, for some choice of $s,p,q$. Here we slightly restrict the function class, insisting that the densities under consideration are bounded and bounded away from 0. In particular, the lower bound condition facilitates the faster minimax estimation rates of Proposition \ref{Prop: minimax estimation rates}; it is shown in \cite{weedEstimationSmoothDensities2019} that removing this condition results in slower rates for most parameter configurations.

\begin{defn}
	Let $1\leq p,q\leq \infty$, $s\geq0$, $B\geq 1$, $M\geq 1\geq m>0$. Define the function class
	\begin{equation}
		\mathcal{F}_{s,p,q}(B;m,M) = \left\{ f\in B^s_{pq}: \int_{\Td} f=1,\quad \|f\|_{B^s_{pq}}\leq B,\quad m\leq f\leq M \,\mathrm{a.e.} \right\};
	\end{equation}
	Note that we always have $1\in\mathcal{F}_{s,p,q}(B;m,M)$, and so the class is non-empty.
	Henceforth we fix $p=2$ and consider $q,B,m,M$ to be given. Define
	$$ \mathcal{F}(s) := \mathcal{F}_{s,2,q}(B;m,M). $$
\end{defn}
For large $s$ and smaller values of $B\geq1$, the condition $f\leq M$ is superfluous. However, the imposition of the uniform lower bound $f\geq m>0$ means that $\mathcal{F}(s)$ is a strict subset of the more typical parameter space $\{f\in B^s_{2q}:f\geq0,\int f = 1,\|f\|_{B^s_{2q}}\leq B\}$. Also, it is clear from the definition (\ref{Eq: Besov norm definition}) that the continuous embedding $B^s_{pq}\subset B^r_{pq}$ holds with operator norm 1, so $\mathcal{F}(s) \subset \mathcal{F}(r)$ for $r\leq s$.

\subsubsection{Notation}

For a probability density $f$, let $P_f$ and $E_f$ denote respectively the probability and expectation when $X_1,\ldots, X_n \overset{\mathrm{i.i.d.}}{\sim}f$. For real numbers $a,b$, we write $a\wedge b=\min(a,b)$ and $a\vee b = \max(a,b)$. Given sequences $(a_n)$ and $(b_n)$, we write $a_n\lesssim b_n$ if there exists a constant $C>0$ that is independent of $n$ such that for all $n$, $a_n \leq Cb_n$; we also write $a_n\simeq b_n$ if $a_n \lesssim b_n$ and $b_n \lesssim a_n$. Given any subset $A$ of a metric space $(\mathcal{A},d)$, we write $|A|_d$ for the $d$-diameter of $A$, defined by
$$ |A|_d := \sup_{x,y\in A} d(x,y).$$
Given a subset $B\subset \mathcal{A}$ and a point $a\in \mathcal{A}$, we define the distance of $a$ to $B$ as
$$ d(a,B) := \inf_{b\in B}d(a,b).$$

\subsection{Description of the Problem}\label{subsection: problem description}

Suppose initially that $f\in\mathcal{F}(r)$ for some given $r\geq0$. We wish to construct a confidence set $C_n$ for the unknown density $f$; informally, we would like $C_n$ to contain $f$ with (some chosen) high probability. Specifically, given $\alpha\in(0,1)$, we require any confidence set $C_n = C_n(\alpha,X_1,\ldots,X_n)$ to have \emph{honest coverage} at level $1-\alpha$ over the class $\mathcal{F}(s)$, that is, there exists $n_0\in\mathbb{N}$ such that for all $n\geq n_0$,
\begin{equation}\label{Eq: honesty definition}
	 \inf_{f\in\mathcal{F}(r)} P_f(f\in C_n) \geq 1-\alpha.
\end{equation}
The `honesty' refers to the uniformity over $\mathcal{F}(r)$. We remark that in the minimax paradigm, one must necessarily insist on honesty, since the true density $f_0$ is unknown: `dishonest' adaptive confidence sets exist (see \cite[Corollary 8.3.10]{gineMathematicalFoundationsInfiniteDimensional2015}), but the index $n_0$ from which coverage is valid depends on the unknown $f$, so such procedures produce questionable guarantees in practice.

It is clear that the smaller the set $C_n$, the more informative it is; otherwise one could just take $C_n$ to be the whole parameter space $\mathcal{F}(r)$. Thus we desire the $W_2$-diameter of our set $C_n$ to shrink as quickly as possible in $n$. Suppose $C_n$ satisfies the honest coverage condition (\ref{Eq: honesty definition}) for some $\alpha\in(0,1)$, and let $r_n$ be a positive sequence such that for some $\beta>\alpha$ and every $n\geq n_0$, we have
\begin{equation}\label{Eq: O_P minimax estimation rate}
	\inf_{\tilde{f}_n}\sup_{f\in\mathcal{F}(r)}P_f\big(W_2(\tilde{f}_n,f)\geq r_n\big)\geq \beta.
\end{equation}
Here, the infimum is taken over all \emph{estimators} (i.e. measurable functions) $\tilde{f}_n = \tilde{f}_n(X_1,\ldots,X_n)$.
Then by Lemma 2 in \cite{robinsAdaptiveNonparametricConfidence2006}, the $W_2$-diameter of $C_n$ satisfies, for $n\geq n_0$,
$$ \sup_{f\in\mathcal{F}(r)}P_f\left(|C_n|_{W_2}\geq r_n\right) \geq \beta - \alpha;$$
in particular, its diameter cannot shrink faster than $r_n$ with high probability. We define the \emph{minimax estimation rate} (in probability) over $\mathcal{F}(s)$, denoted $r_n^*(s)$, to be the `slowest' sequence (i.e. the largest such sequence up to a multiplicative prefactor) $r_n$ such that (\ref{Eq: O_P minimax estimation rate}) is satisfied for some $\beta>0$ and some $n_0\geq1$. Usually this rate depends on the smoothness parameter $s$.

\begin{remark}\label{Remark: expectation vs O_P rates}
	 The term `minimax estimation rate' is often reserved for any sequence $\bar{r}_n$ such that
	$$ \inf_{\tilde{f}_n}\sup_{f\in\mathcal{F}(r)} E_f W_2(\tilde{f}_n,f) \simeq \bar{r}_n. $$
	By Markov's inequality, we have that $r_n^* \lesssim \bar{r}_n$. In fact, as shown by Proposition \ref{Prop: minimax estimation rates} below, in this problem the rates $r_n^*$ and $\bar{r}_n$ coincide (possibly up to a logarithmic factor when $d=2$).
\end{remark}

In general, it is unrealistic to assume that the regularity $r$ is known. Thus we find ourselves in an adaptation problem, where we wish to construct procedures that do not depend on the unknown smoothness $r$, but which result in (near-)optimal performance for a range of values of $r$. In order to highlight the main ideas, let us consider the two class adaptation problem, where for some fixed $s>r\geq0$ we consider the model $\mathcal{F}(r)$, but also seek optimal performance over the smoother subclass $\mathcal{F}(s) \subset \mathcal{F}(r)$. We discuss after Theorem \ref{Thm: d > 2 confidence sets} how one might construct confidence sets adapting to a continuous window of smoothnesses $[r,R]$ or even all $r\geq0$ simultaneously.
\begin{defn}\label{Def: optimal adaptive confidence set}
	We say that $C_n = C_n(\alpha,\alpha', X_1,\ldots,X_n)$ is a \emph{near-optimal adaptive $W_2$ confidence set over $\mathcal{F}(s)\cup\mathcal{F}(r)$}, $s>r$, if it satisfies the following properties, for given $\alpha,\alpha'\in(0,1)$:
	\begin{enumerate}[(i)]
		\item \textbf{Honest Coverage:} for all $n$ sufficiently large,
		\begin{equation}\label{Eq: acs def, coverage}
			\inf_{f\in\mathcal{F}(r)} P_f(f\in C_n) \geq 1-\alpha;
		\end{equation}
		\item \textbf{Diameter Shrinkage:} there exists a constant $K = K(\alpha')>0$ such that
		\begin{equation}\label{Eq: acs def, slow shrinkage}
			\sup_{f\in\mathcal{F}(r)} P_f(|C_n|_{W_2} > KR_n(r)) \leq \alpha'
		\end{equation}
		and
		\begin{equation}\label{Eq: acs def, fast shrinkage}
			\sup_{f\in\mathcal{F}(s)} P_f(|C_n|_{W_2} > KR_n(s)) \leq \alpha' ,
		\end{equation}
		for $n$ large enough, where the rate sequences $R_n(r)$ and $R_n(s)$ satisfy
		$$ R_n(r) \leq a_nr_n^*(r) \quad \text{and} \quad R_n(s) \leq a_nr_n^*(s), $$
		for $r_n^*(r)$ and $r_n^*(s)$ the minimax rates of estimation over $\mathcal{F}(r)$ and $\mathcal{F}(s)$ respectively and $a_n$ some power of $\log{n}$.
	\end{enumerate}
\end{defn}

Typically, for optimal adaptive confidence sets one insists that the rates $R_n(r),R_n(s)$ in (\ref{Eq: acs def, slow shrinkage}) and (\ref{Eq: acs def, fast shrinkage}) are equal up to constants to the minimax estimation rates $r_n^*(r), r_n^*(s)$. Our definition of `near-optimal' allows for $R_n(t)$ to equal $r_n^*(t), t=r,s$, up to a logarithmic factor in $n$, and is thus a slight relaxation. Admitting this relaxation does not alter the (existence and) non-existence results of \cite{bullAdaptiveConfidenceSets2013},  \cite{carpentierHonestAdaptiveConfidence2013}, \cite{hoffmannAdaptiveInferenceConfidence2011}, \cite{gineMathematicalFoundationsInfiniteDimensional2015}, since these results are due to a polynomial discrepancy between minimax estimation and testing rates; see Section \ref{Subsection: nonexistence} below.

We only consider the problem of adaptation in the smoothness parameter and do not address the question of adaptation to other parameters in the definition of the class $\mathcal{F}(s)$, such as the Besov norm bound $B$. See Remark \ref{Remark: adapting over other parameters} below for a discussion of this issue.

\subsection{Adaptive $W_2$ Confidence Sets on $\Td$}

Our first theorem exhaustively classifies the parameter configurations for which adaptive honest confidence sets exist for $W_2$ loss; in the cases where such confidence sets do exist, an explicit construction is given in Theorem \ref{Thm: d > 2 confidence sets} below.

\begin{theorem}\label{Thm: existence and nonexistence of conf sets}
	Fix $1\leq q\leq\infty$, $B\geq1$, $M\geq1\geq m>0$. Consider the two class adaptation problem for confidence sets as defined by (\ref{Eq: acs def, coverage})-(\ref{Eq: acs def, fast shrinkage}).
	\begin{enumerate}[(i)]
		\item Let $d\leq4$ and $s>r\geq0$. Then for any $\alpha,\alpha'>0$, there exists a near-optimal adaptive $W_2$ confidence set.
		\item Let $d>4$ and $0\leq r<s\leq \frac{2d-4}{d-4}r+\frac{d}{d-4}$. Then for any $\alpha,\alpha'>0$, there exists a near-optimal adaptive $W_2$ confidence set.
		\item Let $d>4$ and $0\leq r<s$ with $s> \frac{2d-4}{d-4}r+\frac{d}{d-4}$. Then for any $\alpha,\alpha'>0$ such that $2\alpha + \alpha'<1$, no near-optimal adaptive $W_2$ confidence set exists.
	\end{enumerate}
\end{theorem}
\begin{remark}\label{Remark: choice of W_2}
	We have focussed on the particular choice of $W_2$; by Jensen's inequality, this distance dominates $W_p$ for $1\leq p<2$. Since the minimax estimation rates in these problems are independent of $p$ (c.f. Proposition \ref{Prop: minimax estimation rates}), this means that the above existence results hold for $W_p,1\leq p\leq 2$, in particular for the important case of $W_1$. Moreover, in the case of $W_1$, one may remove the lower bound condition in the definition of $\mathcal{F}(s)$; see Remark \ref{Remark: classes not bounded below for W_1} below.
\end{remark}

Theorem \ref{Thm: existence and nonexistence of conf sets} says that in low dimensions, $d\leq 4$, there exists a confidence set which adapts optimally in $W_2$-diameter to \emph{any} two smoothnesses $s>r\geq0$. As the construction does not depend on $s$, in fact adaptation occurs simultaneously for all $s\geq r$ (strictly speaking, $r\leq s\leq S$ where $S$ is the regularity of the wavelet basis used), where $r$ is a chosen `baseline' smoothness.
Contrast this to the case of $L_p$ loss, $2\leq p\leq \infty$: for $p<\infty$, in any dimension, there exists a (near-)optimal adaptive confidence set if and only if $s\leq \frac{p}{p-1}r$ (\cite{bullAdaptiveConfidenceSets2013}, \cite{carpentierHonestAdaptiveConfidence2013}); for $L_{\infty}$ loss, adaptive confidence sets do not exist for any choice of $s>r\geq0$ (\cite{lowNonparametricConfidenceIntervals1997}, \cite{hoffmannAdaptiveInferenceConfidence2011}). See \cite[Section 8.3]{gineMathematicalFoundationsInfiniteDimensional2015} for a complete account of the $L_2$ and $L_{\infty}$ theory.

In higher dimensions $d>4$, Theorem \ref{Thm: existence and nonexistence of conf sets} gives a `window' of smoothnesses for which adaptation occurs, in a similar vein to the case of $L_p,p<\infty$. However, for the $W_2$ loss the window is significantly wider; moreover, regardless of how small we choose $r\geq0$, this window has width at least $\frac{d}{d-4}$, whereas for $L_p, 2\leq p<\infty$, the window is of width $\frac{r}{p-1}\leq r$, which will be very narrow for small values of $r$.

These results are related to the fact that $W_2$ is a weaker loss function than $L_p$: specifically, Proposition \ref{Prop: W-Besov comparison} and (\ref{Eq: W_2 - log Sobolev comparison}) show that on the class $\mathcal{F}(s)$, $W_2$ is comparable to a Sobolev (or Besov) norm of smoothness -1. In very low dimensions $d=1,2$, the estimation rate is independent of the smoothness parameter $s$, meaning that any confidence set satisfying (\ref{Eq: acs def, slow shrinkage}) automatically satisfies the faster shrinkage condition (\ref{Eq: acs def, fast shrinkage}) (with a possibly enlarged constant $K$). In low dimensions $d=3,4$, one finds a very fast minimax testing separation rate, which can be as fast as the parametric rate of estimation $n^{-1/2}$ (this is implied by the above existence results and Lemma \ref{Lemma: confidence set impossibility} below). Even in higher dimensions, there is a substantial acceleration in the testing separation rate as compared to $L_2$ loss. Meanwhile, although there is also some acceleration in the estimation rates, the effect is not so pronounced. This explains the wider window of adaptation seen in Theorem \ref{Thm: existence and nonexistence of conf sets} for $W_2$ loss, as compared to $L_p$ loss: the greater discrepancy between testing and estimation rates gives more room for adaptation to take place.

Theorem \ref{Thm: existence and nonexistence of conf sets} is proved in Section \ref{Section: Existence Results}; we outline the arguments now. For the existence result, we use the method of constructing confidence sets via risk estimation as in \cite{juditskyNonparametricConfidenceSet2003}, \cite{caiAdaptiveConfidenceBalls2006}, \cite{robinsAdaptiveNonparametricConfidence2006}; see \cite[Section 6.4]{gineMathematicalFoundationsInfiniteDimensional2015} for a concise summary of these ideas. These methods require the loss function under consideration to be a Hilbert space norm. Accordingly, we upper bound $W_2$ by a suitable Sobolev-type norm for which one can perform risk estimation with fast convergence rates; moreover, the estimation rates for this dominating norm differ from those for $W_2$ by only a logarithmic factor. In particular, the notions of near-optimal adaptive confidence sets for these two loss functions are equivalent. The non-existence result is obtained using a testing argument as in \cite{hoffmannAdaptiveInferenceConfidence2011}, \cite{bullAdaptiveConfidenceSets2013} and others, together with a lower bound for the minimax separation rate in a related testing problem. Moreover, the precise characterisation of the separation rate identifies a certain small subset of $\mathcal{F}(r)$ consisting of `problematic' densities which, once removed, permit the existence of confidence sets (with honesty relative to a smaller set of densities), as in the previous two references. We discuss the existence of these more general confidence sets after Theorem \ref{Thm: testing rate}. These theoretical results and constructions extends more generally to the study of adaptive honest confidence sets with negative Sobolev norm distances, and we discuss them in Section \ref{subsection: Sobolev norm theory extension}.
For $p>2$, \cite{carpentierHonestAdaptiveConfidence2013} develops a construction of adaptive $L_p-$confidence sets whose radii are selected via testing. Though an extension of these ideas to $W_p$-confidence sets should be possible, we do not pursue it here as the methodology greatly differs from the one used in the present paper.

\subsection{Adaptive $W_1$ Confidence Sets on $\Rd$}

The case of densities on $\Rd$ is also of great interest; there are several situations in which it is unrealistic to assume compact support of the density $f$. Accordingly, let $X_1,\ldots,X_n$ be an i.i.d. sample drawn from some unknown density $f$ on $\Rd$. We take the Wasserstein-1 distance $W_1$ to be our loss function. We generalise our methods from the case of $\Td$ to produce adaptive confidence sets for $f$ which adapt over similar function classes $\mathcal{G}(s)$, defined in \eqref{Eq: function class definition Rd} below and involving a constant $L$ which uniformly bounds the exponential moments of the densities in $\mathcal{G}(s)$. The discussion following Theorem \ref{Thm: existence and nonexistence of conf sets} is relevant in this context as well: in particular, since the confidence sets constructed in cases (i) and (ii) do not depend on $s$, adaptation in fact takes place for the full range of possible values of $s$ (i.e. $s\geq r$ when $d\leq 4$ and $s$ in some given window when $d>4$). 

\begin{theorem}\label{Thm: existence and nonexistence of conf sets - Rd}
	Fix $1\leq q\leq\infty$, $B\geq1$, $M\geq1\geq m>0$. Consider the two class adaption problem for confidence sets defined by (\ref{Eq: acs def, coverage})-(\ref{Eq: acs def, fast shrinkage}), with function classes $\mathcal{F}$ replaced by $\mathcal{G}$ and $W_2$ in place of $W_1$.
	\begin{enumerate}[(i)]
		\item Let $d\leq4$ and $s>r\geq0$. Then for any $\alpha,\alpha'>0$, there exists a near-optimal adaptive $W_1$ confidence set.
		\item Let $d>4$ and $0\leq r<s\leq \frac{2d-4}{d-4}r+\frac{d}{d-4}$. Then for any $\alpha,\alpha'>0$, there exists a near-optimal adaptive $W_1$ confidence set.
		\item\label{item: last case} Let $d>4$, $L$ be large enough and $0\leq r<s$ with $s> \frac{2d-4}{d-4}r+\frac{d}{d-4}$. Then for any $\alpha,\alpha'>0$ such that $2\alpha + \alpha'<1$, no near-optimal adaptive $W_1$ confidence set exists.
	\end{enumerate}
\end{theorem}

The bound $L$ on exponential moments in \eqref{Eq: function class definition Rd} is a technical condition which allows us to construct adaptive estimators and confidence sets via the method of risk minimization (see Section \ref{Section: Rd}). We are naturally interested in the existence of confidence sets for large $L$, i.e. on larger classes of densities. Moreover, small values of $L$ may lead to empty classes (see the discussion after Definition \ref{def: class on Rd} below) for which the theory of confidence sets is superfluous.

\subsection{Extension to negative Sobolev norm distances}\label{subsection: Sobolev norm theory extension}

To better understand the phenomena in Theorems \ref{Thm: existence and nonexistence of conf sets} and \ref{Thm: existence and nonexistence of conf sets - Rd}, it is elucidating to consider negative order Sobolev norm loss, $H^{-t}=B^{-t}_{22},t>0$ (see Appendix \ref{Section: Wavelet Appendix} for definitions), since the $W_2$ distance is dominated by such a norm (see (\ref{Eq: W_2 - log Sobolev comparison}) below). One finds that the minimax estimation rate for $t\geq d/2$ is (up to a log factor) $n^{-1/2}$, so no meaningful adaptation is required and one constructs a confidence set which `adapts' over all smoothnesses as in Proposition \ref{Prop: confidence set, d=1,2} below. When $t<d/2$, computations analogous to those in Section \ref{Section: Existence Results} show that the gap between testing and estimation rates are wider for larger $t$, enabling adaptation over a larger window of regularities (see Remark \ref{Remark: weak Sobolev norms} below). Here, one finds a continuous transition as $t$ increases from 0 (which is the $L_2$ case) to $d/2$, at which point confidence sets can adapt to any two smoothnesses. However, the specific geometry of the parameter space induced by the loss function is crucial, rather than how weak the loss function is \emph{per se}: if instead we consider $B^{-t}_{\infty\infty}$ loss, when $t<d/2$ the minimax estimation and testing rates can be shown to coincide; meanwhile, the estimation rate is independent of the smoothness parameter when $t\geq d/2$. So in the case of $B^{-t}_{\infty\infty}$ loss, when $t<d/2$ no adaptive confidence sets exist for \emph{any} two smoothnesses by Lemma \ref{Lemma: confidence set impossibility} below, but for $t\geq d/2$ they trivially exist.

Whenever they exist, the construction of confidence sets in Section \ref{Section: Existence Results} below extends easily to the case of negative order Sobolev norms $H^{-t}, t>0$, and other Besov norms using norm embeddings as in \cite[Section 4.3]{gineMathematicalFoundationsInfiniteDimensional2015}; see Remark \ref{Remark: weak Sobolev norms} below.
%The success of the techniques used also suggests that there is insight to be gained from the study of weaker loss functions than the classical choice of $L_p$ norms, replacing these with negative order Sobolev or Besov norms.
% As it transpires from our proof, the statements of Theorem \ref{Thm: existence and nonexistence of conf sets} are also valid, up to minor modifications, if one replaces the Wasserstein distance with one of these weaker norm-distances. For instance, for $t>0$ and the norm of the negative Sobolev space $H_2^{-t}$, full adaptation is possible for dimension $d\leq 4t$, while a restriction to a window of the form $\left(r;(d-4t)^{-1}(dt+(2d-4t)r)\right]$ is necessary for higher dimensions. It is consistent with Theorem \ref{Thm: existence and nonexistence of conf sets} as we compare $W_2$ with norms that have a smoothness parameter $t=1$.
	
	\section{Proof of Theorem \ref{Thm: existence and nonexistence of conf sets}}\label{Section: Existence Results}

\subsection{A Hilbert Norm Upper Bound for $W_2$}

We wish to construct confidence sets by performing risk estimation. The inner product structure of Hilbert space norms makes them particularly amenable to risk estimation, and so we seek some Hilbert norm which upper bounds the $W_2$ distance.

For this, we introduce the \emph{logarithmic Sobolev norm} (\cite[Section 4.4]{gineMathematicalFoundationsInfiniteDimensional2015}; see \cite{castilloNonparametricBernsteinMises2013}, \cite{castilloBernsteinvonMisesPhenomenon2014} for another statistical application of such norms).
\begin{defn}
	Define the $H^{-1,\delta}$ norm of $f\in L_2(\Td)$ as
	$$ \|f\|_{H^{-1,\delta}} = |\langle f,1 \rangle| + \left(\sum_{l\geq0} 2^{-2l}\max(l,1)^{2\delta}\|\langle f, \psi_{l\cdot}\rangle\|_{2}^2\right)^{1/2}.$$
	Note the similarity to the definition of the $B^{-1}_{22} = H^{-1}$ norm given by (\ref{Eq: Besov norm definition}); indeed, when $\delta=0$ the two norms coincide with the Sobolev norm of regularity -1.
	We refer to this as a `logarithmic' Sobolev space because the parameter $\delta$ measures the smoothness of $f$ on a logarithmic scale. 
\end{defn}

We require the following comparison inequality from \cite{weedEstimationSmoothDensities2019}.

\begin{prop}[Theorem 3, \cite{weedEstimationSmoothDensities2019}]\label{Prop: W-Besov comparison}
	Let $1\leq p<\infty$. Let $f,g$ be two densities in $L_p(\Td)$, and assume that for almost every $x\in\Td$, $M\geq\max(f(x), g(x))\geq m>0$, for real numbers $M$ and $m$. Then
	\begin{equation}\label{Eq: W-B estimate}
		M^{-1/p'}\|f-g\|_{B^{-1}_{p\infty}} \lesssim W_p(f,g) \lesssim m^{-1/p'}\|f-g\|_{B^{-1}_{p1}},
	\end{equation} 
	where $\frac{1}{p}+\frac{1}{p'}=1$, and the constants depend only on $d,p$ and the wavelet basis. Moreover, when $p=1$, one may choose $m=0$. 
\end{prop}

This result is an extension of the celebrated Kantorovich-Rubinstein duality formula, which states that for two \emph{probability measures} $\mu,\nu$ on $\Td$,
\begin{equation}\label{Eq: K-R Duality}
	W_1(\mu,\nu) = \sup_{h\in\mathrm{Lip}_1(\Td)} \int h\,\ud(\mu-\nu),
\end{equation}
where the supremum is taken over all functions $h:\Td\to\R$ with Lipschitz constant bounded by 1. We may relate this to (\ref{Eq: W-B estimate}) using the sequence of norm-continuous embeddings (\cite[Section 4.3]{gineMathematicalFoundationsInfiniteDimensional2015})
$$ B^{-1}_{11}\subset \left(B^1_{\infty\infty}\right)^* \subset BL(\Td)^* \subset \left(B^1_{\infty 1}\right)^* \subset B^{-1}_{1\infty},$$
where $BL(\Td)$ is the space of bounded Lipschitz functions on $\Td$ (note that any Lipschitz function on $\Td$ is bounded, so $BL(\Td)$ and $\mathrm{Lip}_1(\Td)$ coincide). However, in order to generalise this to $W_p,p>1$, one must impose that the probability measures have densities which are bounded and bounded away from zero; indeed, for densities not bounded below, no norm provides a similar comparison to $W_p$ (\cite[Theorem 7]{weedEstimationSmoothDensities2019}), and convergence rates are slower than those in Proposition \ref{Prop: minimax estimation rates}. Thus the restriction from the usual choices of Besov norm-balls to the classes $\mathcal{F}(s),s\geq0$ is necessary.

A simple application of the Cauchy-Schwarz inequality confirms that $H^{-1,\delta} \subset B^{-1}_{21}$ as soon as $\delta>1/2$. Thus in conjunction with the upper bound in Proposition \ref{Prop: W-Besov comparison}, we have that, for $r\geq0,$ $f\in\mathcal{F}(r)$ and $\tilde{f}_n$ any estimator of $f$,
\begin{equation}\label{Eq: W_2 - log Sobolev comparison}
	W_2(f,\tilde{f}_n) \lesssim \|f-\tilde{f}_n\|_{B^{-1}_{21}} \lesssim \|f-\tilde{f}_n\|_{H^{-1,\delta}},
\end{equation}
where the first constant depends on the parameters of the class $\mathcal{F}(r)$, but the second constant depends only on the wavelet basis and $d$.

\begin{remark}\label{Remark: classes not bounded below for W_1}
	When using $W_1$ loss, one may consider the class $\mathcal{F}(s)$ with the choice $m=0$, i.e. densities are not required to be bounded away from zero. Then the $H^{-1,\delta}$ norm still provides an upper bound for $W_1$ for densities in $\mathcal{F}(s)$ due to the upper bound in (\ref{Eq: W-B estimate}) and the sequence of continuous embeddings
	$ H^{-1,\delta} \subset B^{-1}_{21} \subset B^{-1}_{11},$
	where the second embedding follows from Jensen's inequality (with operator norm 1). 
\end{remark} 

For the remainder of this section, we work in $H^{-1,\delta}$ risk; as soon as $\delta>1/2$, this provides a Hilbert norm upper bound for the $W_2$ risk. In particular, any coverage guarantee for a $H^{-1,\delta}$ ball is automatically inherited by the $W_2$ ball with the same centre and radius scaled by the embedding constant from (\ref{Eq: W_2 - log Sobolev comparison}). Of course, by constructing confidence sets for a stronger loss function, we may not be able to attain near-optimal diameter shrinkage, but we shall see that this is not the case.

\subsection{Construction of Confidence Sets}

We first give the minimax estimation rates for the problem under consideration. These are important for two reasons: firstly, they provide the benchmark for the `size' of an optimal confidence set. Moreover, our confidence sets are centred at a suitable estimator of $f$, which must perform well for the resulting confidence set to also have good performance. In the density estimation problem, the estimation rates for $W_2$ loss are as follows:
\begin{prop}\label{Prop: minimax estimation rates}
	Let $s\geq0$ and let $r_n^*(s)$ denote the minimax rate of estimation over $\mathcal{F}(s)$. Then
	$$ r_n^*(s) \lesssim \begin{cases}
		n^{-1/2}, \quad &d=1,\\
		n^{-1/2}\log{n}, \quad &d=2,\\
		n^{-\frac{s+1}{2s+d}}, \quad &d\geq3,
	\end{cases}$$
	where the constant depends on the parameters of the class $\mathcal{F}(s)$ and the wavelet basis.
	Moreover, for any $s\geq0$,
	$$ r_n^*(s) \gtrsim \begin{cases}
		n^{-1/2}, \quad &d=1,2\\
		n^{-\frac{s+1}{2s+d}}, \quad &d\geq3,
	\end{cases}$$
	where the infimum is over all estimators $\tilde{f}_n$ based on a sample of size $n$.
\end{prop}
The upper bounds follow from Theorem 1 in \cite{weedEstimationSmoothDensities2019} and Remark \ref{Remark: expectation vs O_P rates}. The lower bounds are proved as in Theorem 6.3.9 in \cite{gineMathematicalFoundationsInfiniteDimensional2015}, where one ensures the existence of a suitable $W_2$-separated set using the lower bound in Proposition \ref{Prop: W-Besov comparison}. See also Theorem 2 in \cite{weedEstimationSmoothDensities2019}.

%Accordingly, for $s\geq0$ define rates 
%\begin{equation}\label{Eq: r_n(s) definition}
%	r_n(s) = \begin{cases}
%		n^{-1/2}, \quad &d=1,2, \\
%		n^{-\frac{s+1}{2s+d}}, \quad &d\geq3.
%	\end{cases}
%\end{equation}
%In the $d=2$ case there is a logarithmic gap in Proposition \ref{Prop: minimax estimation rates}, so we take the lower bound in our definition of $r_n(s)$ (due to our notion of near-optimality, taking the upper bound will not change any of our theorems). However, in all other dimensions, $r_n(s)$ is the precise minimax estimation rate $r_n^*(s)$.

We centre our confidence sets at an estimator $\hat{f}_n$ of $f$ which has near-optimal convergence over the classes $\mathcal{F}(s)$ and $\mathcal{F}(r)$. The theory of adaptive estimation is relatively complete, and in the vast majority of cases it is possible to construct adaptive estimators which converge at the minimax estimation rate (perhaps up to a logarithmic factor) over a wide range of smoothnesses - we mention only the classical references \cite{lepskiiProblemAdaptiveEstimation1991} and \cite{donohoDensityEstimationWavelet1996}. 

The consideration of Wasserstein loss adds a minor complication to the usual case of `norm-type' loss functions. The Wasserstein distance $W_p(f,\tilde{f}_n)$ is only well-defined if $\tilde{f}_n$ is also a density, and thus we ought to insist that any estimator we define is indeed a density almost surely. To achieve this, given any wavelet-based estimator of the form
$$ \tilde{f}_n = \tilde{f}_{-1} + \sum_{l\geq0}\sum_{k=0}^{2^{ld}-1}\tilde{f}_{lk}\psi_{lk} $$
where $\tilde{f}_{lk}$ are the wavelet coefficients of the estimator, we insist that $\tilde{f}_{-1} = 1$. This ensures that $\int_{\Td}\tilde{f}_n = 1$. The problem of non-negativity is more subtle. In \cite{weedEstimationSmoothDensities2019}, it was addressed by projecting $\tilde{f}_n$ onto the class of densities $\mathcal{F}(r)$ with respect to the $B^{-1}_{p1}$ norm, where $r$ is the smallest regularity to which we want to adapt.  Then, for $f\in\mathcal{F}(r)$, denoting the projection of $\tilde{f}_n$ by $\tilde{f}^D_n$,
$$ W_p(f,\tilde{f}^D_n) \lesssim \|f-\tilde{f}^D_n\|_{B^{-1}_{p1}} \leq \|f-\tilde{f}_n\|_{B^{-1}_{p1}} + \inf_{g\in\mathcal{F}(r)}\|\tilde{f}_n-g\|_{B^{-1}_{p1}} \leq 2\|f-\tilde{f}_n\|_{B^{-1}_{p1}},$$
and so it suffices to analyse the performance of $\tilde{f}_n$ in $B^{-1}_{p1}$ loss. However, this projection step makes the estimator essentially intractable. Instead, we use the well-known $L_{\infty}$ consistency of the adaptive estimators considered below (c.f. \cite{donohoDensityEstimationWavelet1996}, for example) together with the fact that the densities in $\mathcal{F}(r)$ are uniformly bounded away from 0 to conclude that for sufficiently large $n$, with high probability $\tilde{f}_n$ is in fact a probability density. Whenever $\tilde{f}_n$ fails to be non-negative, we simply replace it with an arbitrary choice of density (e.g. uniform); as $n\to\infty$, this event occurs with vanishing probability.

\begin{theorem}\label{Thm: Thresholded estimator}
	Let $d\geq2$. Then there exists an estimator $\hat{f}_n$ of $f$ such that for all $n\geq n_0(B)$ and all $s\geq0$,
	$$ \sup_{f\in\mathcal{F}(s)} E_f\|f-\hat{f}_n\|_{H^{-1,\delta}}^2 \lesssim (\log{n})^{2\delta}\left(\frac{n}{\log{n}}\right)^{-\frac{2(s+1)}{2s+d}},$$
	where the constant depends on $B,d$ and the wavelet basis.
\end{theorem}

The definition of $\hat{f}_n$ and proof of Theorem \ref{Thm: Thresholded estimator} can be found in Appendix \ref{Section: additional torus proofs}, and follows from the classical ideas of \cite{donohoDensityEstimationWavelet1996}.

Next, we introduce a $U$-statistic to perform risk estimation.
Recall that given any estimator $\tilde{f}_n$ of $f$ such that $\langle \tilde{f}_n, 1\rangle = 1$, the $H^{-1,\delta}$ loss can be expressed as
$$ \|f-\tilde{f}_n\|_{H^{-1,\delta}}^2 = \sum_{l\geq0}2^{-2l}(l\vee1)^{2\delta}\sum_{k=0}^{2^{ld} -1}\langle f-\tilde{f}_n,\psi_{lk}\rangle^2. $$
To estimate this loss, we use the approach of sample splitting. Suppose we have a sample of size $2n$ which we divide into two subsamples
$$ \Sone=(X_1,\ldots,X_n),\quad \Stwo = (X_{n+1},\ldots,X_{2n}).$$
Denote expectation with respect to sample $i$ by $E^{(i)}$; we denote variances and probabilities accordingly.
We compute our estimator $\tilde{f}_n=\tilde{f}_n(X_1,\ldots,X_n)$ based on $\Sone$ and, for $j\geq0$, define the $U$-statistic based on the sample $\Stwo$ as
\begin{equation}\label{Eq: U-statistic definition}
	U_{n,j}(\tilde{f}_n) = \frac{2}{n(n-1)}\sum_{i<i',i,i'\in\Stwo}\sum_{l<j}2^{-2l}(l\vee1)^{2\delta}\sum_{k=0}^{2^{ld}-1}\left(\psi_{lk}(X_i) - \langle\psi_{lk},\tilde{f}_n\rangle \right)\left(\psi_{lk}(X_{i'}) - \langle\psi_{lk},\tilde{f}_n\rangle\right). 
\end{equation}
Since the sample is i.i.d., we see that
$$ E^{(2)}_fU_{n,j}(\tilde{f}_n) = \sum_{l<j}2^{-2l}(l\vee1)^{2\delta}\sum_{k=0}^{2^{ld}-1}\langle\psi_{lk},f-\tilde{f}_n\rangle^2 = \|K_j(f-\tilde{f}_n)\|_{H^{-1,\delta}}^2.$$
Thus $U_{n,j}(\tilde{f}_n)$ is an unbiased estimator of the $j^{th}$ resolution level approximation of the loss $\|f-\tilde{f}_n\|_{H^{-1,\delta}}$. The key idea behind the $U$-statistic is that the removal of the diagonal in the outermost sum in (\ref{Eq: U-statistic definition}) eliminates the highest variance terms. Thus by averaging over $O(n^2)$ terms with small variance, we expect the $U$-statistic to have very small variance (as in Theorem 6.4.6 of \cite{gineMathematicalFoundationsInfiniteDimensional2015}). This is confirmed by the next lemma.

\begin{lemma}\label{Lemma: U-stat variance bound}
	Assume $f\in L^{\infty}(\Td)$ is a probability density, and $\tilde{f}_n$ is an estimator for $f$ based on the subsample $\Sone$. Then
	\begin{align}
		\mathrm{Var}^{(2)}(U_{n,j}(\tilde{f}_n)) &\leq \frac{4\|f\|_{\infty}}{n}\left(\underset{l\geq -1}{\max}\ 4^{-l} (1\vee l)^{2\delta}\right)\|K_j(f-\tilde{f}_n)\|_{H^{-1,\delta}}^2 + \frac{2\|f\|_{\infty}^2}{n(n-1)}\sum_{l\leq j-1}2^{l(d-4)}(l\vee1)^{4\delta}  \nonumber \\
		&=: \kappa_{n,j,\delta}^2(f). \label{Eq: kappa definition}
	\end{align}
\end{lemma}
This result is analogous to Theorem 4.1 in \cite{robinsAdaptiveNonparametricConfidence2006}; for completeness, we give a proof in Appendix \ref{Section: additional torus proofs}.

With the adaptive estimator $\hat{f}_n$ and the $U$-statistic $U_{n,j}(\hat{f}_n)$ in hand, we are now ready to give the construction of optimal confidence sets for the two-class adaptation problem. 

We first note that for $d=1,2$, the minimax rates of estimation from Proposition \ref{Prop: minimax estimation rates} do not depend on the smoothness parameter $s$; in particular, the two diameter shrinkage conditions (\ref{Eq: acs def, slow shrinkage}) and (\ref{Eq: acs def, fast shrinkage}) become a single condition. Thus in these dimensions, defining an adaptive confidence set is very easy; indeed, there is no meaningful adaptation which needs to take place.

When $d=1$, the empirical measure is a minimax optimal estimator of the sampling measure (see, for instance, \cite{weedSharpAsymptoticFinitesample2019} or \cite{fournierRateConvergenceWasserstein2015}). When $d=2$, we centre at the adaptive estimator from Theorem \ref{Thm: Thresholded estimator} in place of the empirical measure $P_n$, as $P_n$ is no longer minimax optimal, and standard kernel or wavelet projection estimators require choices of tuning parameters depending on the smoothness parameter to attain optimal rates.
\begin{prop}\label{Prop: confidence set, d=1,2}
	\begin{enumerate}[(i)]
	\item Let $d=1$. Consider the two-class adaptation problem over $\mathcal{F}(s)\cup\mathcal{F}(r)$ where $s>r\geq0,q\in[1,\infty],B\geq1,M\geq1\geq m>0$ are all fixed. Then given any $\alpha\in(0,1)$, the confidence set based on a sample $X_1,\ldots,X_n$ defined by
	$$ C_n = \left\{ g\in\mathcal{F}(r) : W_2(P_g,P_n) \leq D\alpha^{-1/2}n^{-1/2} \right\} $$
	is an optimal adaptive $W_2$ confidence set, where $P_n = n^{-1}\sum_{i=1}^n \delta_{X_i}$ is the $n$-sample empirical measure and the constant $D$ depends on $B, m$ and the wavelet basis. 
	\item Let $d=2$. Consider the two-class adaptation problem over $\mathcal{F}(s)\cup\mathcal{F}(r)$ where $s>r\geq0,q\in[1,\infty],B\geq1,M\geq1\geq m>0$ are all fixed. Then given any $\alpha\in(0,1)$, the confidence set based on a sample $X_1,\ldots,X_n$ defined by
	$$ C_n = \left\{ g\in\mathcal{F}(r) : W_2(g,\hat{f}_n) \leq D\alpha^{-1/2}n^{-1/2}(\log{n})^{2+\delta} \right\} $$
	is a near-optimal adaptive $W_2$ confidence set, where $\hat{f}_n$ is the adaptive estimator from Theorem \ref{Thm: Thresholded estimator} and the constant $D$ depends on $B, m$ and the wavelet basis.
	\end{enumerate}
\end{prop}
The diameter shrinkage conditions are met trivially, while honest coverage follows from Chebyshev's inequality in a standard fashion.

When $d\geq3$, the minimax rates depend on the smoothness parameter and so the diameter shrinkage condition differs between $\mathcal{F}(r)$ and $\mathcal{F}(s)$, $r\neq s$. In particular, this precludes any confidence set $C_n$ with deterministic radius, as used above. Instead, we centre at the adaptive estimator $\hat{f}_n$ from Theorem \ref{Thm: Thresholded estimator}, and use the estimate of its loss provided by the $U$-statistic $U_{j,n}(\hat{f}_n)$ as defined in (\ref{Eq: U-statistic definition}) to determine the radius. We write $U_{j} := U_{j,n}(\hat{f}_n)$ in the sequel.

\begin{theorem}\label{Thm: d > 2 confidence sets}
	Let $d\geq3$. Fix $B\geq1,M\geq1\geq m>0,1\leq q\leq\infty$, and let $s>r\geq0$. If $d>4$, assume additionally that $s\leq \frac{2d-4}{d-4}r + \frac{d}{d-4}$. Fix $\alpha\in(0,1)$, and $\delta>1/2$. Consider the confidence set based on a sample of size $2n$, $\Sone\cup\Stwo$ given by
	\begin{equation}\label{Eq: Wasserstein confidence set definition}
		C_n = \left\{ g\in\mathcal{F}(r): \|g-\hat{f}^T_n\|_{H^{-1,\delta}} \leq \sqrt{z_{\alpha}\kappa_{n,j_n,\delta}(g) + U_{j_n} + G(j_n)} \right\}
	\end{equation}
	where $\hat{f}^T_n$ is computed on $\Sone$, $U_{j_n}$ is computed on $\Stwo$ and:
	\begin{itemize}
		\item $ \kappa_{n,j,\delta}^2(g) := \frac{4\|g\|_{\infty}}{n}\|K_j(g-\hat{f}^T_n)\|_{H^{-1,\delta}}^2 + \frac{2\|g\|_{\infty}^2}{n(n-1)}\sum_{l\leq j-1}2^{l(d-4)}(l\vee1)^{4\delta}$;
		\item $j_n$ is such that $2^{j_n}\simeq \left(\frac{n}{\log{n}}\right)^{\frac{1}{2r+d/2}}$;
		\item $G(j_n)$ = $j_n^{2\delta}2^{-2j_n(r+1)}\log{n}$;
		\item $z_{\alpha}$ = $(\alpha/2)^{-1/2}$.
	\end{itemize}
	Then for all $n\geq n_0(B)$, $C_n$ satisfies (\ref{Eq: acs def, coverage}), as well as (\ref{Eq: acs def, slow shrinkage}) and (\ref{Eq: acs def, fast shrinkage}) for a suitable constant $K>0$ depending on $r,s,\alpha,\alpha'$ and the parameters of the class $\mathcal{F}(r)$ with the rates
	$$ R_n(r) = (\log{n})^{\delta+\frac{r+1}{2r+d}}n^{-\frac{r+1}{2r+d}}, \qquad R_n(s) = (\log{n})^{\delta+\frac{s+1}{2s+d}}n^{-\frac{s+1}{2s+d}}. $$
	In particular, $C_n$ is a near-optimal adaptive $W_2$ confidence set over $\mathcal{F}(s)\cup\mathcal{F}(r)$.
	\begin{remark}[Adaptation over ranges of classes]
		Note that the construction of $C_n$ is completely independent of $s$, and $\hat{f}_n$ adapts simultaneously over all $s\geq0$. So when $d\leq4$, $C_n$ adapts simultaneously over all $s\geq r$, and when $d>4$, $C_n$ adapts simultaneously over the full window of admissible values of $s$.
	\end{remark}
	\begin{remark}[Adaptation to other parameters]\label{Remark: adapting over other parameters}
		We note that the construction of the confidence set in Theorem \ref{Thm: d > 2 confidence sets} does not depend on $B$ or $m$, and so in fact this particular confidence set is also adaptive over $B\geq1$ and $m>0$, in the sense that any dependence of the minimax rates $r_n^*(r),r_n^*(s)$ on $B$ or $m$ are eventually accounted for by the logarithmic term in $R_n(r),R_n(s)$. (Note however that the constants in our theoretical guarantees explode as $B\to\infty$ or $m\to0$.) However, the construction of $C_n$ does depend on $M$. See \cite{bullAdaptiveConfidenceSets2013} for more discussion on the role of $M$.
	\end{remark}
	\begin{remark}[Adapting to wider ranges of smoothnesses in high dimensions]
		In the $d>4$ case, following the ideas in \cite{bullAdaptiveConfidenceSets2013}, one may still obtain adaptation over a window of the form $[0,R]$ for arbitrary $R>0$ at the cost of removing certain troublesome portions of the classes $\mathcal{F}(r),r\in[0,R]$. In this restricted model, one can identify the smoothness of the unknown density within a window of the form $\left[r,\frac{2d-4}{d-4}r+\frac{d}{d-4}\right]$ using tests as in \cite{bullAdaptiveConfidenceSets2013} or \cite{nicklSharpAdaptiveConfidence2016}. Once this window is identified, in particular the relevant value of $r$, one can use the associated confidence set as constructed in Theorem \ref{Thm: d > 2 confidence sets}.
	\end{remark}
	\begin{remark}(Necessity of log-factors)
		One may ask whether it is possible to remove the log-factors in the shrinkage rates and construct a confidence set with $R_n(r)=r_n^*(r), R_n(s)=r_n^*(s)$. These log factors fundamentally arise from the use of the embedding $H^{-1,\delta}\hookrightarrow B^{-1}_{21}$ for $\delta>1/2$. For confidence sets constructed via risk estimation we conjecture that this is a necessary step, as it is precisely the accelerated risk estimation for Hilbert space norms which enables the adaptivity of the confidence set. Working with the $B^{-1}_{21}$ norm directly, which is in some sense an $L_1$-type norm, it seems one will run into problems as outlined in \cite{lepskiEstimationLrNorm1999} and \cite{caiAccuracyAssessmentHighDimensional2016}, where it is shown that for the $L_1$ norm, risk estimation cannot be performed (polynomially) more accurately than the size of the risk itself in both a Gaussian regression model and a sparse high-dimensional linear model. While we do not yet have any precise negative results for the $B^{-1}_{21}$ norm, risk estimation is itself an important topic of study and thus this question should be addressed in the future. However, it is conceivable that another approach, such as the testing method of \cite{carpentierHonestAdaptiveConfidence2013}, could be used to construct $W_2$ confidence sets with sharp diameter shrinkage rates.
	\end{remark}
	\begin{remark}[Weak Sobolev norms $H^{-t},t>0$]\label{Remark: weak Sobolev norms}
		Our methods extend to the use of negative order Sobolev norms $H^{-t}=B^{-t}_{22}, t>0$ as loss functions in place of $H^{-1,\delta}$ (see Appendix \ref{Section: Wavelet Appendix} for definitions). The analysis of the estimator $\hat{f}_n$ is completely analogous, and one must suitably augment the $U$-statistic $U_{n,j}$ to estimate the $H^{-t}$ loss. One finds that the resulting confidence set $\tilde{C}_n$ adapts to any two smoothnesses $0\leq r<s<\infty$ when $t\geq d/4$; if instead $t<d/4$, adaptation is possible over a window of smoothnesses $0\leq r<s\leq \frac{d}{d-4t}t + \frac{2d-4t}{d-4t}r$. Moreover, in this latter case, the arguments of Section \ref{Subsection: nonexistence} below can be augmented to show that if $s$ does not lie in this window, then no such confidence set can exist.
	\end{remark}
\end{theorem}
The proof of this theorem proceeds similarly to that of Proposition 2.1 in \cite{robinsAdaptiveNonparametricConfidence2006}, and is given in Appendix \ref{Section: additional torus proofs}.

The confidence sets constructed above prove statements (i) and (ii) of Theorem \ref{Thm: existence and nonexistence of conf sets}. 

\subsection{Testing rates and non-existence of Confidence Sets}\label{Subsection: nonexistence}

We turn now to proving the impossibility result (iii) in Theorem \ref{Thm: existence and nonexistence of conf sets}.

The question of existence of adaptive confidence sets is closely related to a composite hypothesis testing problem. This connection was identified in the first works on adaptive confidence sets; for a complete decision-theoretic overview, see \cite[Chapter 8]{gineMathematicalFoundationsInfiniteDimensional2015}. For $\rho\geq0$ and $s>r\geq0$, define the separated function class
$$ \tilde{\mathcal{F}}(r, \rho) := \left\{ f\in\mathcal{F}(r): W_2(f,\mathcal{F}(s)) \geq \rho \right\}$$
We may have $\rho=0$, in which case $\tilde{\mathcal{F}}(r,0) = \mathcal{F}(r)$. However, if $\rho>0$ then $\tilde{\mathcal{F}}(r,\rho)$ is a strict subset of $\mathcal{F}(r)$, disjoint from $\mathcal{F}(s)$. The testing problem we consider is
\begin{equation}\label{Eq: testing problem}
	H_0: f\in\mathcal{F}(s) \quad \mathrm{vs.} \quad H_1: f\in\tilde{\mathcal{F}}(r,\rho).
\end{equation}
As the usefulness of a test is naturally assessed by the sum of its Type I and Type II errors, the minimax rate of testing for the problem \eqref{Eq: testing problem} is defined as any sequence $\left(\rho_n^*\right)_{n\geq 1}$ such that
\begin{itemize}
	\item For any $\beta'>0$, there exists a constant $L=L\left(\beta'\right)$ and a measurable test $\Psi_n: \left(\mathbb{T}^d\right)^n\to \{0,1\}$ such that
	\begin{equation}\label{Eq: uniform consistency}
		\underset{f\in \mathcal{F}(s)}{\sup} \mathbb{E}_f \left[\Psi_n\right] +  \underset{f\in \tilde{\mathcal{F}}(r,L\rho^*_n)}{\sup} \mathbb{E}_f \left[1-\Psi_n\right] \leq \beta'.
	\end{equation}
	\item There exists some $\beta>0$ such that for all $\rho_n=o\left(\rho_n^*\right),$
	\begin{equation}\label{lower condition tests}
		\underset{n\to\infty}{\lim\inf} \ \underset{\Psi_n}{\inf}\Big[\underset{f\in\mathcal{F}(s)}{\sup} \mathbb{E}_f \left[\Psi_n\right] +  \underset{f\in \tilde{\mathcal{F}}(r,\rho_n)}{\sup} \mathbb{E}_f \left[1-\Psi_n\right] \Big] \geq \beta,
	\end{equation}
	where the infimum ranges over the set of tests $\Psi_n$.
\end{itemize}

The following result characterises the role of the minimax testing rate $\rho_n^*$ in the existence and non-existence of confidence sets. Essentially, it says $\rho_n^*$ provides a `speed limit' on how quickly the confidence set can shrink when $f$ is in the smoother submodel $\mathcal{F}(s)$:
\begin{lemma}[Proposition 8.3.6, \cite{gineMathematicalFoundationsInfiniteDimensional2015}] \label{Lemma: confidence set impossibility}
	Let $\rho_n^*$ be the minimax testing rate for (\ref{Eq: testing problem}), and $\Tilde{r}_n(s), \Tilde{r}_n(r)$ be two sequences such that $\Tilde{r}_n(s)=o\left(\rho_n^*\right)$ and $\Tilde{r}_n(s)=o\left(\Tilde{r}_n(r)\right)$. Let $\alpha, \alpha'>0$. Then, for any $\rho_n=o\left(\rho_n^*\right)$ and $L>0$, there does not exist any set $C_n\left(\alpha, X_1,\dots,X_n\right)$ satisfying
	\begin{itemize}
		\item $\liminf_{n\to\infty} \inf_{f\in\mathcal{F}(s)\cup\tilde{\mathcal{F}}(r,\rho_n)} P_f(f\in C_n) \geq 1-\alpha,$
		\item $ \limsup_{n\to\infty}  \sup_{f\in\tilde{\mathcal{F}}(r,\rho_n )}P_f\big(|C_n|_{W_2}>L\Tilde{r}_n(r)\big) \leq \alpha',$
		\item $ \limsup_{n\to\infty}  \sup_{f\in\mathcal{F}(s)}P_f\big(|C_n|_{W_2}>L\Tilde{r}_n(s)\big) \leq \alpha',$
	\end{itemize}
	as long as $\alpha,\alpha'$ are such that $0<2\alpha+\alpha'<\beta$, with $\beta$ as in \eqref{lower condition tests}.
\end{lemma}
This non-existence phenomenon occurs because any $C_n$ satisfying the conditions of the Lemma induces a test
$$ \Psi_n = \ind\{ C_n \cap \tilde{\mathcal{F}}(r,\rho'_n)\neq \varnothing \} $$
which is uniformly consistent for the separation rate $\rho_n'$ in the sense of (\ref{Eq: uniform consistency}) whenever $\rho_n = o(\rho'_n)$. If we were able to choose $\rho'_n$ to be $o(\rho^*_n)$, this would contradict the definition of the minimax testing rate $\rho^*_n$; thus no such confidence set can exist. Note that the argument works for any rate $\tilde{r}_n(s) = o(\rho^*_n)$, not just the minimax rate of estimation; in particular, we can multiply the minimax estimation rate by a poly-logarithmic factor so long as there is a polynomial gap between the testing and estimation rates.

It remains to determine the minimax rate of testing for the problem (\ref{Eq: testing problem}); this is done in the following theorem.
\begin{theorem}\label{Thm: testing rate}
	Assume $s>r\geq0$ and $d> 4$. Let $\rho_n^*$ be the minimax rate of testing for the problem (\ref{Eq: testing problem}). Then there exist a constant $c>0$ depending on the parameters of the class $\mathcal{F}(s)$ and the wavelet basis, and $n_0 = n_0(B,M)$ such that for all $n\geq n_0$, 
	$$\rho_n^* \geq c n^{-\frac{r+1}{2r+d/2}}.$$
	Also, \eqref{lower condition tests} holds for any $\beta<1$.
\end{theorem}

The proof of Theorem \ref{Thm: testing rate} is given in Appendix \ref{Section: additional torus proofs}, and follows a multiple-testing lower bound.
Assume now that $d>4$ and $s>\frac{2d-4}{d-4}r + \frac{d}{d-4}$. Then the minimax rate of testing $\rho^*_n$ is slower than the minimax estimation rate $r_n^*(s)$ by a polynomial factor; in light of Lemma \ref{Lemma: confidence set impossibility}, this means there is no near-optimal adaptive $W_2$ confidence set over $\mathcal{F}(s)\cup\mathcal{F}(r)$ for any practical choice of $\alpha,\alpha'$ (for such a set to exist, we would require $2\alpha+\alpha'\geq1$). This proves statement (iii) of Theorem \ref{Thm: existence and nonexistence of conf sets}. However, this does not rule out the existence of confidence sets satisfying weaker conditions than those in Definition \ref{Def: optimal adaptive confidence set}, namely those listed in Lemma \ref{Lemma: confidence set impossibility} for some $\rho_n\geq L \rho_n^*$, $L>0$. Such sets actually exists in view of Proposition 8.3.7 of \cite{gineMathematicalFoundationsInfiniteDimensional2015} and Theorem \ref{Thm: Thresholded estimator}.

Moreover, the confidence set $C_n$ constructed in Theorem \ref{Thm: d > 2 confidence sets} in conjunction with the argument used to prove Lemma \ref{Lemma: confidence set impossibility} shows that the lower bound of Theorem \ref{Thm: testing rate} is sharp up to a poly-logarithmic factor.
	
	\section{Extension of the Theory to $\mathbb{R}^d$}\label{Section: Rd} 

Having provided a fairly complete resolution of the problem of adaptive $W_2$ confidence sets when the sample space is $\Td$, we extend our results to the case of the unbounded sample space $\Rd$ with $W_1$ loss. The key tool is the Kantorovich-Rubinstein duality formula (\cite{kantorovich1958space})
\begin{equation}\label{Eq: K-R duality, 2}
	W_1(f,g) = \sup_{h\in\mathrm{Lip}_1(\Rd)} \int_{\Rd} h(x)(f(x) - g(x))\,\ud x,
\end{equation}
where $\mathrm{Lip}_1(\Rd)$ is the set of 1-Lipschitz functions on $\Rd$.

We also briefly discuss what happens when using $W_p$ loss for $p>1$ in Section \ref{Subsection: W_p,p>1 on Rd}; the unbounded sample space introduces complications which preclude a direct generalisation of the ideas from Section \ref{Section: Existence Results}.

In this section, it is assumed that we observe $X_1,\ldots,X_n\overset{\mathrm{i.i.d.}}{\sim}f_0$ for some density $f_0$ on $\Rd$, and we wish to perform inference on $f_0$ using $W_1$ as the loss function.

\subsection{Parameter Spaces}

We use an $S$-regular tensor product wavelet basis of $L^2(\Rd)$ of the form
$$ \left\{\phi_k,\psi_{lk}: k\in\Zd,l\geq0\right\}$$
as introduced in Appendix \ref{Section: Wavelet Appendix} (we index the $\psi_{lk}$ using only $k,l$ by a slight abuse of notation). We write $K_j(f)$ for the projection of $f$ onto the first $j$ resolution layers, as in (\ref{Eq: j level projection definition}). Besov norms on $\Rd$, also defined in Appendix \ref{Section: Wavelet Appendix}, are defined analogously to those on $\Td$, and the relation (\ref{Eq: Besov wavelet coefficient bound}) holds.

Our goal is to construct an adaptive confidence set for the true density $f_0$ using risk estimation, where the adaptation occurs with respect to the smoothness parameter $s$. We shall consider functions in $B^s_{2q}$. Unlike our previous classes $\mathcal{F}(s)$ on $\Td$, we need not assume that our densities are bounded away from zero, or something analogous such as sufficiently slow decay in the tails. However, in order to deal with the unboundedness of the sample space $\Rd$, we must impose a moment condition.

For $\alpha,\beta>0$, define the $\alpha,\beta$-exponential moment of a density $f$ as
\begin{equation}\label{Eq: exponential moment definition}
	\mathcal{E}_{\alpha,\beta}(f) := \int_{\Rd} \exp{(\beta\|x\|^{\alpha})}f(x)\,\ud x = E_f\left(e^{\beta \|X\|^{\alpha}}\right).
\end{equation}

\begin{defn}\label{def: class on Rd}
	Let $1\leq p,q\leq\infty,s\geq0,B\geq1, M>0,\alpha,\beta>0$ and $L\geq1$. Define the function class
	\begin{equation}\label{Eq: function class definition Rd}
		\mathcal{G}_{s,p,q}(B,M;\alpha,\beta,L) = \left\{f\in B^{s}_{pq}(\Rd): \int_{\Rd}f =1,\quad \|f\|_{B^s_{pq}}\leq B,\quad 0\leq f\leq M \,\text{a.e.},\quad \mathcal{E}_{\alpha,\beta}(f) \leq L \right\}.
	\end{equation}
	Henceforth, we fix $p=2$ and consider $q,B,M,\alpha,\beta,L$ to be given. Define
	$$ \mathcal{G}(s) := \mathcal{G}_{s,2,q}(B,M;\alpha,\beta,L). $$
\end{defn}
Observe that for $M$ close to 0 and $L$ close to 1, the class $\mathcal{G}(s)$ is empty. We therefore assume in the sequel that $L$ is sufficiently large (depending on $M,B$) for $\mathcal{G}(s)$ to be non-empty.

The focus on $p=2$ is quite natural in view of the material developed in the previous section, relying on risk estimation to compute the diameter of confidence sets. Combining the exponential moment condition and the bound on the $B^s_{2q}$-norm, we prove in Lemma \ref{lemma: class inclusion} that densities in $\mathcal{G}(s)$ also have their $B^s_{1q}$-norm bounded by a constant depending on the class parameters.

\subsection{Estimation Upper Bounds for $W_1$}\label{Subsection: Rd estimation}

As before, we should insist on our estimator $\tilde{f}_n$ being a density almost surely. Indeed, the fact that $\tilde{f}_n$ has total mass 1 is vital to the proof of Proposition \ref{Prop: W1 loss decomposition} below. However, we note that there is no intrinsic requirement in (\ref{Eq: K-R duality, 2}) that $f$ and $g$ should be nonnegative, and so we will allow our estimators to take negative values. If a genuine density is required, one can just take the positive part of the estimator and renormalize. 

The following proposition gives an upper bound on the $W_1$ distance which is convenient for wavelet estimators.
\begin{prop}\label{Prop: W1 loss decomposition}
	For any probability density $f$ with a finite first moment and any estimator $\tilde{f}_n$ of $f$ which has a finite first moment almost surely, we have that
	\begin{equation}\label{Eq: W1 loss decomposition}
		W_1(\tilde{f}_n,f) \lesssim \sum_{k\in\Zd} \|k\||\langle f-\tilde{f}_n,\phi_k \rangle| + \sum_{l\geq0} 2^{-l\left(\frac{d}{2} + 1\right)} \sum_{k\in\Zd} |\langle f-\tilde{f}_n, \psi_{lk}\rangle|,
	\end{equation}
	where the constant depends only on the wavelet basis.
\end{prop}
\begin{remark}
	Let $\hat{f}_n$ be some estimator of $f$, not necessarily with total mass 1. We obtain an estimator which integrates to 1 almost surely, which we call $\tilde{f}_n$, by renormalising the first wavelet layer of $\hat{f}_n$, that is, renormalising $\hat{f}_0 := K_0(\hat{f}_n)$. Then we set
	$$ \tilde{f}_n = \frac{\hat{f}_0}{\int 
		\hat{f}_0(x)\,\ud x} + \sum_{l\geq0}\sum_{k\in\Zd}\langle \hat{f}_n,\psi_{lk}\rangle \psi_{lk}.$$
	Note that while one can perform this procedure for any estimator $\hat{f}_n$, it is particularly simple for wavelet-based estimators. Assuming $L_1$-consistency of $\hat{f}_n$, $\hat{f}_0\to K_0(f)$ and thus $K_0(\tilde{f}_n)\to K_0(f)$ in $L_1$. Moreover, for the wavelet estimators we use below, this convergence occurs very fast, at the rate $n^{-\frac{S}{2S+d}}$, where $S$ is the regularity of the wavelet basis. Thus it suffices to consider the un-normalised estimator $\hat{f}_n$ in the decomposition (\ref{Eq: W1 loss decomposition}) whenever $s\leq S-1$, which we do in the sequel.
\end{remark} 

We first establish an upper bound for the estimation rate over the class $\mathcal{G}(s)$.
\begin{theorem}\label{Thm: Rd estimation upper bound}
	For any $s\geq0$, there exists an estimator $\hat{f}_n$ such that for all sufficiently large $n$,
	$$ \sup_{f\in\mathcal{G}(s)} E_f\,W_1(\hat{f}_n,f) \lesssim \begin{cases}
		(\log{n})^{\frac{\gamma d}{2}+1}n^{-1/2}, \quad &d=2, \\
		(\log{n})^{\frac{\gamma d}{2}}n^{-\frac{s+1}{2s+d}}, \quad &d\geq3.
	\end{cases} $$
	where $\gamma$ is a constant depending on $\alpha$ and $\beta$ only, and the constant depends on the parameters of the class $\mathcal{G}(s)$ and the wavelet basis. For $d=1$, the empirical measure $P_n$ satisifies \[  \sup_{f\in\mathcal{G}(s)} E_f\,W_1(P_n,P_f) \lesssim n^{-1/2}.\]
\end{theorem}
\begin{remark}
	These rates are sharp up to a logarithmic factor so long as $L$ is sufficiently large: one uses a reduction to a multiple testing problem as in the proof of the lower bounds in Proposition \ref{Prop: minimax estimation rates}, and then uses an analogous collection of well-separated densities defined on some common compact set. For large enough $L$, the compact support ensures that these densities have suitable exponential moments and so belong to $\mathcal{G}(s)$.
\end{remark}
\begin{remark}
	An inspection of the proof reveals that in fact it suffices to assume a suitable polynomial moment, depending on $s$; however, for convenience we assume an exponential moment which works for all $s\geq0$.
\end{remark}
The proofs of Proposition \ref{Prop: W1 loss decomposition} and Theorem \ref{Thm: Rd estimation upper bound} are given in Appendix \ref{Section: Rd proofs}. The estimator $\hat{f}_n$ is simply a wavelet projection estimator which is zero outside of a growing compact set; the risk outside of the compact is controlled using the moment assumption.

As in the case of $\Td$, we require an adaptive estimator.
\begin{theorem}\label{Thm: Rd adaptive estimator}
	Let $d\geq2$, and let $\gamma>0$ be as in Theorem \ref{Thm: Rd estimation upper bound}. Then there exists an estimator $\hat{f}_n$ of $f$ such that for all $n\geq n_0(B)$ and all $s\geq0$,
	$$ \sup_{f\in\mathcal{G}(s)} E_f\, W_1(\hat{f}_n,f) \lesssim (\log{n})^{\frac{\gamma d}{2}}\left(\frac{n}{\log{n}}\right)^{-\frac{s+1}{2s+d}}, $$
	where the constant depends on the parameters of the class $\mathcal{G}(s)$ and the wavelet basis.
\end{theorem}
The definition of $\hat{f}_n$ and proof of Theorem \ref{Thm: Rd adaptive estimator} are given in Appendix \ref{Section: Rd proofs}.

\subsection{Construction of Confidence Sets}

Let us now concretely state the two-class adaptation problem we wish to solve. Fix two smoothnesses $s>r\geq0$ and consider the model $\mathcal{G}(r) = \mathcal{G}(r)\cup\mathcal{G}(s)$. Given $\alpha\in(0,1)$, we seek a confidence set $C_n$ which has honest coverage at level $1-\alpha$, that is, for all $n$ sufficiently large,
\begin{equation}\label{Eq: Rd honest coverage}
	\inf_{f\in\mathcal{G}(r)}P_f(f\in C_n) \geq 1-\alpha,
\end{equation}
as well as the two diameter shrinkage conditions: for all $\alpha'>0$ there exists a constant $K=K(\alpha')>0$ such that
\begin{align}
	\sup_{f\in\mathcal{G}(r)} P_f(|C_n|_{W_1}>KR_n(r)) &\leq \alpha', \label{Eq: Rd diameter shrinkage slow} \\
	\sup_{f\in\mathcal{G}(s)} P_f(|C_n|_{W_1}>KR_n(s)) &\leq \alpha', \label{Eq: Rd diameter shrinkage fast}
\end{align}
where $R_n(r)$ and $R_n(s)$ equal the convergence rates in Theorem \ref{Thm: Rd estimation upper bound} up to a poly-logarithmic factor.

As discussed previously, the $d=1$ and $d=2$ cases are straightforward given the existence of the estimator from Theorem \ref{Thm: Rd adaptive estimator}, since here the convergence rates do not depend on the smoothness $r$. We thus restrict our attention to the case $d\geq3$.

Let $X_1,\ldots,X_{2n}$ be an i.i.d. sample from the unknown $f\in\mathcal{G}(r)$. We split the sample as before into two equal halves, $\Sone$ and $\Stwo$, and denote by $P^{(i)},E^{(i)}$ probabilities and expectations taken over $\mathcal{S}_i$. We wish to construct a confidence set via risk estimation, centred at the estimator $\hat{f}_n$ which we compute using $\Sone$. Proposition \ref{Prop: W1 loss decomposition} provides a natural upper bound for $W_1(f,\hat{f}_n)^2$ which we then decompose into several terms. Define the thresholds $\kappa_{-1n} = \kappa_{0n} \simeq (\log{n})^{\gamma}, \kappa_{ln} = 2^l\kappa_{0n}$ for $\gamma$ chosen as in Theorem \ref{Thm: Rd estimation upper bound}. Applying the Cauchy-Schwarz inequality several times, we obtain the bound
\begin{align}\label{Eq: Rd square loss decomposition}
	W_1(f,\hat{f}_n)^2 \leq& 3\Bigg( (\log{n})^{\gamma(d+2)}\left[\sum_{\|k\|_{\infty}\leq \kappa_{-1n}}\langle f-\hat{f}_n,\phi_k\rangle^2 + j\sum_{l<j}2^{-2l}\sum_{\|k\|_{\infty}\leq\kappa_{ln}}\langle f-\hat{f}_n,\psi_{lk}\rangle^2 \right] \nonumber \\
	&\ldots + \left[\sum_{l\geq j}2^{-l\left(\frac{d}{2}+1\right)}\sum_{\|k\|_{\infty}\leq\kappa_{ln}}|\langle f-\hat{f}_n,\psi_{lk}\rangle|\right]^2 \nonumber \\
	&\ldots + \left[\sum_{\|k\|_{\infty}> \kappa_{-1n}}\|k\||\langle f,\phi_k\rangle| + \sum_{l\geq0}2^{-l\left(\frac{d}{2}+1\right)}\sum_{\|k\|_{\infty}>\kappa_{ln}}|\langle f, \psi_{lk}\rangle| \right]^2 \Bigg).
\end{align}
The final term is controlled using the moment assumption on $f\in\mathcal{G}(r)$; indeed, from the proof of Theorem \ref{Thm: Rd estimation upper bound} we have that for all $f\in\mathcal{G}(r)$, this term is bounded above by
\begin{equation}\label{Eq: Rd remainder term definition}
	\Delta_{n} := C(d)L^2(\log{n})^{2\gamma}n^{-1},
\end{equation}
where $C(d)$ is a constant depending only on $d$ and the wavelet basis.

We next consider the remaining terms in (\ref{Eq: Rd square loss decomposition}). We introduce pseudo-distances $\tilde{W}^{(n,j)}(f,g)$ defined as
\begin{equation}\label{Eq: W_1 proxy distances definition}
	\tilde{W}^{(n,j)}(f,g) = \left[\sum_{\|k\|_{\infty}\leq \kappa_{-1n}}\langle f-g,\phi_k\rangle^2 + j\sum_{l<j}2^{-2l}\sum_{\|k\|_{\infty}\leq\kappa_{ln}}\langle f-g,\psi_{lk}\rangle^2 \right]^{1/2} + \sum_{l\geq j}2^{-l\left(\frac{d}{2}+1\right)}\sum_{\|k\|_{\infty}\leq\kappa_{ln}}|\langle f-g,\psi_{lk}\rangle|.
\end{equation}
Observe that for $f,g\in\mathcal{G}(r)$,
$$ W_1(f,g) \leq \sqrt{3(\log{n})^{\gamma(d+2)}}\cdot \tilde{W}^{(n,j)}(f,g) + \sqrt{3\Delta_n};$$
this is true uniformly over $r\geq0$. Since $\sqrt{\Delta_n}$ converges (up to a logarithmic factor) at the parametric rate, this means that any diameter shrinkage condition with respect to $\tilde{W}^{(n,j)}$ provides an analogous shrinkage condition for $W_1$, with only a slightly worse rate. Moreover, the first part of $\tilde{W}^{(n,j)}(f,g)$ is well-suited to estimation using a $U$-statistic. To this end, define the $U$-statistic
\begin{align} \label{Eq: Rd U-stat definition}
	V_{n,j} = V_{n,j}(\hat{f}_n) :=& \frac{2}{n(n-1)}\sum_{i<i',i,i'\in\Stwo} \Bigg[ \sum_{\|k\|_{\infty}\leq \kappa_{-1n}}\left(\phi_k(X_i) - \langle \hat{f}_n,\phi_k\rangle \right)\left(\phi_k(X_{i'}) - \langle \hat{f}_n,\phi_k\rangle \right) \nonumber \\ &\ldots+ j\sum_{l<j}2^{-2l}\sum_{\|k\|_{\infty}\leq\kappa_{ln}}\left(\psi_{lk}(X_i) - \langle \hat{f}_n,\psi_{lk}\rangle \right)\left(\psi_{lk}(X_{i'}) - \langle \hat{f}_n,\psi_{lk}\rangle \right) \Bigg].
\end{align}
Clearly we have that $E^{(2)}_f V_{n,j}$ is equal to the square of the first term in (\ref{Eq: W_1 proxy distances definition}). Analogously to Lemma \ref{Lemma: U-stat variance bound}, one shows that $V_{n,j}$ has small variance.
\begin{lemma}\label{Lemma: Rd U-stat variance bound}
	For $f\in L_{\infty}(\Rd)$, we have that, for some constant $C_d$ depending only on $d$ and the wavelet basis,
	\begin{align*}
		\mathrm{Var}^{(2)}_f(V_{n,j}) &\leq \frac{C_d}{2}\left(\frac{j^2\|f\|_{\infty}^2(\log{n})^{\gamma d}}{n(n-1)}\sum_{l<j}2^{l(d-4)} + \frac{\|f\|_{\infty}}{n}\left[\sum_{\|k\|_{\infty}\leq \kappa_{-1n}}\langle f-\hat{f}_n,\phi_k\rangle^2 + j^2\sum_{l<j}2^{-4l}\sum_{\|k\|_{\infty}\leq\kappa_{ln}}\langle f-\hat{f}_n,\psi_{lk}\rangle^2\right]\right) \\
		&\leq C_d\left(\frac{j^2\|f\|_{\infty}^2(\log{n})^{\gamma d}}{n(n-1)}\sum_{l<j}2^{l(d-4)} + \tilde{W}^{(n,j)}(f,\hat{f}_n)^2 \right) \\
		&=: \lambda_{j,n}^2(f).
	\end{align*}
\end{lemma}

For the second part of $\tilde{W}^{(n,j)}(f,\hat{f}_n)$, we use the concentration arguments from the proof of Theorem \ref{Thm: Rd adaptive estimator} to show that this term is suitably small with high probability uniformly over $f\in\mathcal{G}(r)$.

Given a sequence $(j_n)$, we write $\tilde{W}^{(n)}$ for $\tilde{W}^{(n,j_n)}$, and $V_{j_n}$ for $V_{n,j_n}$.

\begin{theorem}\label{Thm: Rd adaptive confidence sets}
	Let $d\geq3$. Fix $B\geq1,M>0,\alpha,\beta,L>0,1\leq q\leq \infty,$ and $s>r\geq0$. Let $\gamma\geq1$ be as in Theorem \ref{Thm: Rd estimation upper bound}. If $d>4$, assume additionally that $s\leq\frac{2d-4}{d-4}r + \frac{d}{d-4}$. Fix $\alpha\in(0,1)$. Consider the confidence set based on a sample of size $2n$ given by
	\begin{equation}\label{Eq: Rd conf set definition}
		C_n = \left\{g\in\mathcal{G}(r): \tilde{W}^{(n)}(g,\hat{f}_n) \leq C(d)\sqrt{z_{\alpha}\lambda_{n,j_n}(g)+V_{j_n}+ G_{j_n}}\right\}
	\end{equation}
	where $\hat{f}_n$ is computed on $\Sone$, $V_{j_n}$ is computed on $\Stwo$, $C(d)$ is a constant depending on $d$ and the wavelet basis, and:
	\begin{itemize}
		\item $\lambda_{n,j_n}(g)$ is as in Lemma \ref{Lemma: Rd U-stat variance bound};
		\item $j_n$ is such that $2^{j_n}\simeq\left(\frac{n}{\log{n}}\right)^{\frac{1}{2r+d/2}}$;
		\item $G_{j_n} = (\log{n})^{\gamma d + 1}2^{-2j_n(r+1)}$;
		\item $z_{\alpha} = (\alpha/2)^{-1/2}$.
	\end{itemize}
	Then for all $n\geq n_0(B)$, $C_n$ satisfies (\ref{Eq: Rd honest coverage}), as well as (\ref{Eq: Rd diameter shrinkage slow}) and (\ref{Eq: Rd diameter shrinkage fast}) for a suitable constant $K>0$ with the rates
	$$ R_n(r) = (\log{n})^{\gamma(d+1)}\left(\frac{n}{\log{n}}\right)^{-\frac{r+1}{2r+d}}, \quad R_n(s) = (\log{n})^{\gamma(d+1)}\left(\frac{n}{\log{n}}\right)^{-\frac{s+1}{2s+d}}. $$
	In particular, $C_n$ is a near-optimal adaptive $W_1$ confidence set over $\mathcal{F}(s)\cup\mathcal{F}(r)$.
\end{theorem}
The proof is almost identical to that of Theorem \ref{Thm: d > 2 confidence sets}; a more detailed argument can be found in Appendix \ref{Section: Rd proofs}. In particular, this proves statements (i) and (ii) of Theorem \ref{Thm: existence and nonexistence of conf sets - Rd}.

\subsection{Non-Existence of Confidence Sets}

We now turn to the non-existence result (iii) in Theorem \ref{Thm: existence and nonexistence of conf sets - Rd}, a consequence of Lemma \ref{Lemma: confidence set impossibility} (which holds in a general decision theoretic framework). We therefore require a lower bound on the minimax separation rate in the testing problem
\begin{equation}\label{Eq: testing problem Rd}
	H_0: f\in\mathcal{G}(s) \quad \mathrm{vs.} \quad H_1: f\in\tilde{\mathcal{G}}(r,\rho),
\end{equation}
where the separated alternative $\tilde{\mathcal{G}}(r,\rho)$ is defined analogously to before:
$$ \tilde{\mathcal{G}}(r,\rho) := \left\{f\in\mathcal{G}(r): W_1(f,\mathcal{G}(s))\geq\rho\right\}. $$
\begin{theorem}\label{Thm: testing rate Rd}
	Assume that $d>4$ and $s>r\geq0$. Let $\rho_n^*$ be the minimax rate of testing for the problem (\ref{Eq: testing problem Rd}). Then, for $L$ sufficiently large in (\ref{Eq: function class definition Rd}), there exist a constant $c>0$ depending on the parameters of the class $\mathcal{G}(s)$ and the wavelet basis, and $n_0 = n_0(B,M)$ such that for all $n\geq n_0$, 
	$$\rho_n^* \geq c n^{-\frac{r+1}{2r+d/2}}.$$
	Also, \eqref{lower condition tests} holds for any $\beta<1$.\\
\end{theorem}

The proof is given Appendix \ref{Section: Rd proofs}, and is similar to the proof of Theorem \ref{Thm: testing rate}. As before, this implies statement (iii) of Theorem \ref{Thm: existence and nonexistence of conf sets - Rd}.

\subsection{The Case of $W_p,p>1$}\label{Subsection: W_p,p>1 on Rd}

We briefly explain why the above techniques do not extend to other Wasserstein distances $W_p$, for $p>1$.

On $\Rd$, one may bound the $W_1$-distance using the Kantorovich-Rubinstein duality formula (\ref{Eq: K-R Duality}). However, the generalisation of this formula on $\Td$ for $W_p,p>1$, Proposition \ref{Prop: W-Besov comparison}, relies crucially on the compactness of $\Td$ and the fact that the densities to which it applies are bounded uniformly away from zero, say by $m>0$; moreover, the constant in the upper bound is inversely proportional to a power of $m$. To apply this result, certainly we would have to consider only positive densities, use Proposition \ref{Prop: W-Besov comparison} on some compact, and then use a moment condition to control terms outside of this compact.
However, a stronger moment condition will lead to a smaller lower bound $m$ over large compacts; this antagonistic relationship cannot be resolved without a polynomial contribution to the convergence rates.
Given that the $W_1$ estimation rates on $\Rd$ in Theorem \ref{Thm: Rd estimation upper bound} match the rates from the compact case up to logarithmic factors, as well as the fact that such rates for $W_p$ on $\Td$ do not depend on $p$, we conjecture that this method does not lead to sharp upper bounds.

	\subsection*{Acknowledgments}
	
	 The authors gratefully thank Ismaël Castillo and Richard Nickl for their guidance in this project and their careful reading of this paper.
	
	\appendix
	
	\section{Proofs for Section \ref{Section: Existence Results}}\label{Section: additional torus proofs}

We first give the definition of our adaptive estimator. The estimator is based on the empirical wavelet coefficients, defined as
$$ \hat{f}_{lk} := \frac{1}{n}\sum_{i=1}^n \psi_{lk}(X_i).$$
We also write $f_{l\cdot}$ and $\hat{f}_{l\cdot}$ for the vectors of coefficients $(f_{lk}:0\leq k<2^{ld})$ and $(\hat{f}_{lk}:0\leq k<2^{ld})$ respectively.

Next, define the truncation point $l_{\max}$ such that
$$ 2^{l_{\max}} \simeq \left(\frac{n}{\log{n}}\right)^{1/d}, $$
and for $0\leq l\leq l_{\max}$, define the thresholds
$$ \tau_{l} := \tau 2^{\frac{ld}{2}}\left(\frac{\log{n}}{n}\right)^{1/2}, $$
for some $\tau>0$ to be chosen below, depending only on $B,d,M$ and the wavelet basis. We then define
\begin{equation}\label{Eq: thresholded estimator definition}
	\hat{f}_n := 1 + \sum_{l=0}^{l_{\max}} \ind_{\left\{\|\hat{f}_{l\cdot}\|_2^2 > \tau_{l}^2\right\}} \sum_{k=0}^{2^{ld}-1}\hat{f}_{lk}\psi_{lk}.
\end{equation}

To prove Theorem \ref{Thm: Thresholded estimator}, we must first collect some results on the expectation and concentration of the empirical wavelet coefficients $\hat{f}_{lk}$.
\begin{lemma}\label{Lemma: H-J inequality}
	Let $f\in\mathcal{F}(s)$ and let $\hat{f}_{lk}$ be the empirical wavelet coefficients of $f$ based on a sample of $n$ observations. Then for every $t\geq2$ there exists a constant $C_t$ depending only on $t$ such that for all $l\geq0$ satisfying $2^{ld}\leq n$,
	$$ E\left| \hat{f}_{lk} - f_{lk}\right|^t \leq C_t M\|\psi\|_{\infty}^{t-2}n^{-t/2}. $$
\end{lemma}
For $t=2$, the proof is immediate from the i.i.d. assumption on the data, the orthonormality of the wavelets and the bound $\|f\|_{\infty}\leq M$. For $t>2$, the result follows from the $t=2$ case and Hoffmann-J\o rgensen's inequality (\cite{hoffmann-jorgensenSumsIndependentBanach1974}, \cite[Theorem 3.1.22]{gineMathematicalFoundationsInfiniteDimensional2015}).
We also require a concentration result for the $\hat{f}_{lk}$; for this we use Bernstein's inequality (\cite[Theorem 3.1.7]{gineMathematicalFoundationsInfiniteDimensional2015}).
\begin{prop}\label{Prop: Bernstein's inequality}[Bernstein's Inequality]
	Let $Y_1,\ldots,Y_n$ be independent centred random variables which are almost surely bounded by $c>0$ in absolute value. Let $\sigma^2 = n^{-1}\sum_{i=1}^nEY_i^2$ and $S_n = \sum_{i=1}^n Y_i$. Then for all $u\geq0$,
	$$ P(|S_n|>u) \leq 2\exp{\left(-\frac{u^2}{2n\sigma^2 + \frac{2cu}{3}}\right)}. $$
\end{prop}
For fixed $l,k$ and $f\in\mathcal{F}(s)$, the random variables $(\psi_{lk}(X_i) - f_{lk})$ are i.i.d., centred, bounded by $2^{ld/2}\|\psi\|_{\infty} =: c_l$, and have variance bounded by $M$. Thus from  Bernstein's inequality, we deduce that
\begin{equation}\label{Eq: concentration for wavelet coeffs}
	P_f\left(|\hat{f}_{lk}-f_{lk}|>u\right) \leq 2\exp\left(-\frac{nu^2}{2M + \frac{2c_lu}{3}}\right).
\end{equation}

We also need a result on wavelet approximations in the $H^{-1,\delta}$ norm to control bias terms. The following lemma about the error of $j$-level approximations to Besov functions is standard; see Propositions 4.3.8 and 4.3.14 in \cite{gineMathematicalFoundationsInfiniteDimensional2015}, for instance.
\begin{lemma}\label{Lemma: Besov projection accuracy}
	Let $0\leq s<S$ and $1\leq q\leq \infty$, $\delta\in\R$. Then for $f\in B^{s}_{2q}$, we have that
	\begin{equation}\label{Eq: Besov projection accuracy}
		\|K_j(f) - f\|_{H^{-1,\delta}} \leq C \underset{l\geq j}{\sup}\left(2^{-l(s+1)}l^{\delta}\right) \|f\|_{B^{s}_{2q}},
	\end{equation}
	where the constant $C$ depends only on the wavelet basis. In particular, for $j\geq1\vee \frac{\delta}{s+1}$, we have that
	$$ \|K_j(f) - f\|_{H^{-1,\delta}} \leq C 2^{-j(s+1)}j^{\delta}\|f\|_{B^{s}_{2q}} $$
\end{lemma}

\begin{proof}[Proof of Theorem \ref{Thm: Thresholded estimator}]
	Fix $f\in\mathcal{F}(s)$. Define $l_n(s)$ such that
	$$ 2^{l_n(s)} \simeq B^{\frac{1}{s}}\left(\frac{n}{\log{n}}\right)^{\frac{1}{2s+d}};$$
	for all sufficiently large $n$ depending on $B$, we have that $l_n(s)<l_{\max}$. We then decompose the risk as follows:
	\begin{align}
		\|f-\hat{f}_n\|_{H^{-1,\delta}}^2 =& \sum_{l=0}^{l_n(s)}2^{-2l}(l\vee1)^{2\delta}\|\langle f-\hat{f}_n,\psi_{l\cdot}\rangle\|_2^2 + \sum_{l=l_n(s)+1}^{l_{\max}}2^{-2l}l^{2\delta}\|\langle f-\hat{f}_n,\psi_{l\cdot}\rangle\|_2^2 \nonumber \\
		&+ \sum_{l>l_{\max}}2^{-2l}l^{2\delta}\|\langle f,\psi_{l\cdot}\rangle\|_2^2 \nonumber \\
		=:&\,\, I + II + III. \label{Eq: thresholded estimator proof, decomposition}
	\end{align}
	This is a bias-stochastic decomposition, where we have further divided the stochastic term into terms $I$ and $II$.
	
	We first deal with the bias term $III$: a direct application of Lemma \ref{Lemma: Besov projection accuracy} gives
	\begin{align*}
		III &= \|K_{l_{\max}}(f) - f\|_{H^{-1,\delta}}^2 \\
		&\lesssim l_{\max}^{2\delta}2^{-2l_{\max}(s+1)} \\
		&= o\left((\log{n})^{2\delta}\left(\frac{n}{\log{n}}\right)^{-\frac{2(s+1)}{2s+d}}\right)
	\end{align*}
	for a constant depending on $B$ and the wavelet basis.
	
	Next, we deal with term $I$. For any $l\geq0$, by the triangle inequality we have that
	$$ \|\langle f-\hat{f}_n,\psi_{l\cdot}\rangle \|_2 \leq \|f_{l\cdot} - \hat{f}_{l\cdot}\|_2 + \|\hat{f}_{l\cdot}\|_2\ind_{\left\{ \|\hat{f}_{l\cdot}\|_2\leq \tau_l \right\}} \leq \|f_{l\cdot} - \hat{f}_{l\cdot}\|_{2} + \tau2^{ld/2}\sqrt{\frac{\log{n}}{n}}.$$
	Using Lemma \ref{Lemma: H-J inequality} to control the expectation of the square of the first term, we see that
	\begin{align*}
		E_f(I) &\lesssim \sum_{l=0}^{l_n(s)} 2^{-2l}(l\vee1)^{2\delta}\left[2^{ld}n^{-1} + \tau^22^{ld}\frac{\log{n}}{n}\right] \\
		&\lesssim \tau^2\frac{\log{n}}{n}(l_n(s))^{2\delta}\sum_{l=0}^{l_n(s)}2^{l(d-2)},
	\end{align*}
	for $n$ large enough. Note that $l_n(s)\lesssim \log{n}$. Thus when $d=2$, the sum contributes at most some power of $\log{n}$, and so $E_f(I)$ is clearly sufficiently small. For $d>2$, the final term dominates the sum and so using the definition of $l_n(s)$,
	$$	E_f(I) \lesssim \tau^2(\log{n})^{2\delta}\left(\frac{n}{\log{n}}\right)^{-\frac{2(s+1)}{2s+d}} $$
	as required.
	
	Lastly, we must analyse term $II$. Since we consider resolution levels $l>l_n(s)$, we have that
	$$ \|f_{l\cdot}\|_2 \leq B2^{-ls} < B2^{-l_n(s)s} \simeq \left(\frac{n}{\log{n}}\right)^{-\frac{s}{2s+d}},$$
	for a constant depending only on $B$. Moreover,
	$$ \tau_l = \tau2^{ld/2}\left(\frac{n}{\log{n}}\right)^{-1/2} > \tau2^{l_n(s)d/2}\left(\frac{n}{\log{n}}\right)^{-1/2} \geq \tau \left(\frac{n}{\log{n}}\right)^{-\frac{s}{2s+d}}, $$
	and so for $\tau$ chosen sufficiently large depending only on $B$, we have that $\|f_{l\cdot}\|_2 \leq \tau_l/2$.
	Define events
	$$ A_{l,n} := \left\{ \|\hat{f}_{l\cdot}\|_2 \leq \tau_l \right\}, \quad l_n(s)<l\leq l_{\max}. $$
	Then by the above observations, the triangle inequality, a union bound and the bound (\ref{Eq: concentration for wavelet coeffs}), we have that
	\begin{align}
		P_f(A_{l,n}^c) &\leq P_f(\|\hat{f}_{l\cdot}-f_{l\cdot}\|_2>\tau_l/2) \nonumber \\
		&\leq \sum_{k=0}^{2^{ld}-1}P_f\left(|\hat{f}_{lk} - f_{lk}|>\frac{\tau}{2}\sqrt{\frac{\log{n}}{n}}\right) \nonumber \\
		&\leq 2^{ld}\cdot 2\exp\left(-\frac{\tau^2n\log{n}/4}{2Mn + \tau c_l\sqrt{n\log n}/3}\right) \nonumber \\
		&\lesssim \frac{n}{\log{n}}\exp\left(-C\tau\log{n}\right), \label{Eq: threshold low prob event}
	\end{align}
	for $\tau$ large enough depending on $M$ and the wavelet basis, as $l\leq l_{\max}$ and so $2^{l}\leq (n/\log{n})^{1/d}$. Here, $C$ is an absolute constant. Note that on the event $A_{l,n}^c$, $\langle \hat{f}_n, \psi_{lk}\rangle = \hat{f}_{lk}$, whereas on $A_{l,n}$, $\langle \hat{f}_n,\psi_{lk}\rangle = 0$. Thus for $l_n(s)<l\leq l_{\max}$,
	\begin{align}
		E_f\|\langle \hat{f}_n - f,\psi_{l\cdot}\rangle\|_2^2\ind_{A_{l,n}} &\leq \|\langle f,\psi_{l,\cdot}\rangle\|_2^2 \lesssim 2^{-2ls} \label{Eq: term II first part}
	\end{align}
	for some constant depending on $B$, using (\ref{Eq: Besov wavelet coefficient bound}). Next, using Cauchy-Schwarz in conjunction with (\ref{Eq: threshold low prob event}) and Lemma \ref{Lemma: H-J inequality},
	\begin{align}
		E_f\|\langle \hat{f}^T_n - f,\psi_{l\cdot}\rangle \|_2^2\ind_{A_{l,n}^c} &= \sum_{k=0}^{2^{ld}-1}E_f|\hat{f}_{lk} - f_{lk}|^2\ind_{A_{l,n}^c} \nonumber \\
		&\leq \sum_{k=0}^{2^{ld}-1} \left(E_f|\hat{f}_{lk}-f_{lk}|^4\right)^{1/2}\left(P_f(A_{l,n}^c)\right)^{1/2} \nonumber \\
		&\lesssim 2^{ld}(n\log{n})^{-1/2}n^{-C\tau/2}. \label{Eq: term II second part}
	\end{align}
	Combining the estimates (\ref{Eq: term II first part}) and (\ref{Eq: term II second part}), we may bound $II$ as follows:
	\begin{align*}
		E_f(II) &\lesssim \sum_{l=l_n(s)+1}^{l_{\max}} 2^{-2l}l^{2\delta}\left[2^{-2ls} + 2^{ld}(n\log{n})^{-1/2}n^{-C\tau/2}\right] \\
		&\lesssim (\log{n})^{2\delta}\left[2^{-2(s+1)l_n(s)} + (n\log{n})^{-1/2}n^{-C\tau/2}\sum_{l\leq l_{\max}}2^{l(d-2)} \right].
	\end{align*}
	By the definition of $l_n(s)$, the first term is of the correct order. It remains to consider the second term. When $d=2$ the sum contributes a logarithmic factor and so the second term is clearly sufficiently small. When $d>2$, the sum is dominated by its final term and so the second term inside the brackets is of order
	$$ (n\log{n})^{-1/2}n^{-C\tau/2}2^{l_{\max}(d-2)} \simeq \log{n}^{-1/2}n^{\frac{1}{2} - \frac{2}{d} - \frac{C\tau}{2}};$$
	by choosing $\tau$ sufficiently large, we can make this term sufficiently small for all $s\geq0$. This concludes the proof.
\end{proof}

We will also later require the following lemma, which gives control of the $B^s_{2q}$ norm of the estimator $\hat{f}_n$.
\begin{lemma}\label{Lemma: norm control of adaptive estimator}
	Under the hypotheses of Theorem \ref{Thm: Thresholded estimator}, given $\alpha\in(0,1)$ there exists $n_0 = n_0(\alpha)$ such that for all $n\geq n_0$ and any $f\in\mathcal{F}(s)$, with $P_f$-probability at least $1-\alpha$,
	$$ \|\hat{f}_n\|_{B^s_{2q}} \lesssim B + \tau B^{d/2s}, $$
	where the constant depends on $d,q$ only.
	\begin{proof}
		Let $l_n(s),A_{l,n}$ be as in the previous proof. Further define events $B_{l,n} = \{\|\hat{f}_{l\cdot}-f_{l\cdot}\|_2 \leq \tau_l\}$, and 
		$$A_n = \left(\bigcap_{0\leq l\leq l_n(s)}B_{l,n} \right)\bigcap \left(\bigcap_{l_n(s)<l\leq l_{\max}}A_{l,n}\right).$$
		We have from (\ref{Eq: threshold low prob event}), which holds with $B_{l,n}$ in place of $A_{l,n}$ when $l\leq l_n(s)$, and a union bound that
		$$ P_f(A_n^c) \lesssim l_{\max} \frac{n}{\log{n}}\exp\left(-C\tau\log{n}\right) \lesssim n \exp\left(-C\tau\log{n}\right) $$
		and so by choosing $\tau>0$ sufficiently large (independently of $\alpha$), we can make this smaller than $\alpha$ for all sufficiently large $n$. Then on the event $A_n$, using $(a+b)^q \leq 2^{q-1}(a^q + b^q)$,
		\begin{align*}
			\|\hat{f}_n\|_{B^s_{2q}}^q &= 1 + \sum_{l=0}^{l_{\max}}2^{lqs}\ind_{\left\{\|\hat{f}_{l\cdot}\|_2 > \tau_{l}\right\}}\|\hat{f}_{l\cdot}\|_2^q \\
			&\lesssim 1 + \sum_{l=0}^{l_n(s)}2^{lqs}\|f_{l\cdot}\|_2^q + \sum_{l=0}^{l_n(s)}2^{lqs}\|\hat{f}_{l\cdot} - f_{l\cdot}\|_2^q \\
			&\leq \|f\|_{B^s_{2q}}^q + \sum_{l=0}^{l_n(s)}2^{lqs}\tau_l^q \\
			&= B^q + \tau^q\left(\frac{\log{n}}{n}\right)^{q/2}\sum_{l=0}^{l_n(s)}2^{lq\left(\frac{d}{2}+s\right)} \\
			&\lesssim B^q + \tau^qB^{dq/2s},
		\end{align*}
		by choice of $l_n(s)$, since the sum is dominated by its largest term.
	\end{proof}
\end{lemma}

\begin{proof}[Proof of Lemma \ref{Lemma: U-stat variance bound}]
	The kernel of the $U$-statistic is
	$$ R(x,y) = \sum_{l\leq j-1}2^{-2l}(l\vee1)^{2\delta}\sum_{k=0}^{2^{ld}-1}\left[(\psi_{lk}(x) - \langle \psi_{lk},\tilde{f}_n \rangle)(\psi_{lk}(y) - \langle \psi_{lk}, \tilde{f}_n \rangle)\right]$$
	which is symmetric, and so has Hoeffding decomposition (see Section 11.4 of \cite{vaartAsymptoticStatistics2000})
	\begin{equation}\label{Eq: decomposition parts}
		\begin{aligned}
			U_n(\tilde{f}_n) - E^{(2)}_f U_n(\tilde{f}_n) &= \frac{2}{n}\sum_{i\in\Stwo}(\pi_1R)(X_i) + \frac{2}{n(n-1)}\sum_{i<i',i,i'\in\Stwo}(\pi_2R)(X_i,X_{i'}) \\
			&=: L_n + D_n,
		\end{aligned}
	\end{equation}
	with linear kernel
	$$ (\pi_1R)(x) = \sum_{l\leq j-1}2^{-2l}(l\vee1)^{2\delta}\sum_{k=0}^{2^{ld}-1}\left[(\psi_{lk}(x) - \langle \psi_{lk},f\rangle)\langle \psi_{lk},f-\tilde{f}_n\rangle\right] $$
	and degenerate kernel
	$$ (\pi_2R)(x,y) = \sum_{l\leq j-1}2^{-2l}(l\vee1)^{2\delta}\sum_{k=0}^{2^{ld}-1}\left[ (\psi_{lk}(x) - \langle \psi_{lk},f\rangle)(\psi_{lk}(y) - \langle \psi_{lk},f\rangle) \right]. $$
	One checks that $L_n$ and $D_n$ are uncorrelated. It thus remains to bound their variances separately. For $\mathrm{Var}^{(2)}(L_n)$, we use the uncentred version of the kernel $\pi_1R$ and orthonormality of the wavelet basis
	\begin{align*}
		\mathrm{Var}^{(2)}(L_n) &\leq \frac{4}{n}\int\left(\sum_{l\leq j-1}2^{-2l}(l\vee1)^{2\delta}\sum_{k=0}^{2^{ld}-1}\psi_{lk}(x)\langle\psi_{lk},f-\tilde{f}_n\rangle\right)^2 f(x)\,\ud x \\
		&\leq \frac{4\|f\|_{\infty}}{n}\left(\underset{l\geq -1}{\max}\ 4^{-l} (1\vee l)^{2\delta}\right)\sum_{l\leq j-1}2^{-2l}(l\vee1)^{2\delta}\sum_{k=0}^{2^{ld}-1}\langle \psi_{lk}, f-\tilde{f}_n\rangle^2 \\
		&= \frac{4\|f\|_{\infty}}{n} \left(\underset{l\geq -1}{\max}\ 4^{-l} (1\vee l)^{2\delta}\right) \|K_j(f-\tilde{f}_n)\|_{H^{-1,\delta}}^2.
	\end{align*}
	We next bound $\mathrm{Var}^{(2)}(D_n)$. By the degeneracy of the kernel, the summands are uncorrelated. So
	\begin{align*}
		\mathrm{Var}^{(2)}(D_n) &\leq E^{(2)}\left( \frac{2}{n(n-1)}\sum_{i<i',i,i'\in\Stwo}(\pi_2R)(X_i,X_{i'}) \right)^2 \\
		&\leq \frac{2}{n(n-1)}E^{(2)}_f\left(\sum_{l\leq j-1}2^{-2l}(l\vee1)^{2\delta}\sum_{k=0}^{2^{ld}-1}[\psi_{lk}(X_i)\psi_{lk}(X_{i'})]\right)^2 \\
		&\leq \frac{2\|f\|_{\infty}^2}{n(n-1)} \sum_{l\leq j-1} 2^{-4l}(l\vee1)^{4\delta}\sum_{k=0}^{2^{ld}-1}\left( \int \psi_{lk}(x)^2\,\ud x\right)^2 \\
		&= \frac{2\|f\|_{\infty}^2}{n(n-1)}\sum_{l\leq j-1} 2^{l(d-4)}(l\vee1)^{4\delta},
	\end{align*}
	using the orthonormality of the wavelet basis. Combining these two estimates concludes the proof.
\end{proof}

\begin{proof}[Proof of Theorem \ref{Thm: d > 2 confidence sets}]
	We first establish the coverage condition (\ref{Eq: acs def, coverage}). By Lemma \ref{Lemma: norm control of adaptive estimator}, for all $n$ sufficiently large we have with $P_f$-probability at least $1-\alpha/2$ that $\hat{f}_n$ is in a $B^s_{2q}$-norm ball of constant radius. Thus for any $f\in\mathcal{F}(r)$, with $P_f$-probability at least $1-\alpha/2$, for $n\geq n_0(B,\alpha)$ we have from (\ref{Eq: Besov projection accuracy}) that
	$$\| K_{j_n}(f-\hat{f}_n) - (f-\hat{f}_n)\|_{H^{-1,\delta}}^2 \leq G(j_n). $$
	By conditioning on this event, we have that
	\begin{align*}
		P_f(f\in C_n) &= P_f \left(U_{n,j}(\hat{f}_n) - \|f-\hat{f}_n\|_{H^{-1,\delta}}^2 \geq -G(j) - z_{\alpha}\kappa_{n,j,\delta}(f)\right) \\
		&\geq \left(1-\frac{\alpha}{2}\right)P^{(2)}_f\left(U_{n,j}(\hat{f}_n) - \|K_j(f-\hat{f}_n)\|_{H^{-1,\delta}}^2 \geq -z_{\alpha}\kappa_{n,j,\delta}(f)\right) \\
		&\geq \left(1-\frac{\alpha}{2}\right)\left(1 - \frac{\mathrm{Var}^{(2)}_f(U_{n,j}(\hat{f}_n))}{(z_{\alpha}\kappa_{n,j,\delta}(f))^2}\right) \\
		&\geq \left(1-\frac{\alpha}{2}\right)^2 \\
		&\geq 1-\alpha
	\end{align*}
	by Chebyshev's inequality and Lemma \ref{Lemma: U-stat variance bound}.
	
	We now move on to checking the diameter shrinkage conditions (\ref{Eq: acs def, slow shrinkage}) and (\ref{Eq: acs def, fast shrinkage}). Writing $S_j:= \sum_{l<j}2^{l(d-4)}(l\vee1)^{4\delta}$ and using the fact that for positive numbers $a,b$, $\sqrt{a+b}\leq\sqrt{a}+\sqrt{b}$, for $g\in\mathcal{F}(r)$ we have that
	$\kappa_{n,j_n,\delta}(g) \leq 2\sqrt{M}n^{-1/2}\|g-\hat{f}_n\|_{H^{-1,\delta}} + 2M\sqrt{S_{j_n}}n^{-1}$
	and so $g\in C_n$ if and only if
	$$ \|g-\hat{f}_n\|_{H^{-1,\delta}} \leq \sqrt{z_{\alpha}\frac{2M}{n}\sqrt{S_{j_n}} + U_{j_n}+G(j_n)} + n^{-1/4} \sqrt{2z_{\alpha}\sqrt{M}}\sqrt{\|g-\hat{f}_n\|_{H^{-1,\delta}}}.$$
	For positive numbers $x,a,b$, the inequality $x\leq b+a\sqrt{x}$ implies that $x\leq 2b+2a^2$. Thus the diameter of $C_n$ is bounded by a multiple of
	$$ n^{-1/2}S_{j_n}^{1/4} + \sqrt{U_{j_n}} + \sqrt{G_{j_n}} + n^{-1/2}.$$
	We consider each of these terms separately; note that the final term is always sufficiently small.
	
	First, consider $G(j_n)$: this is deterministic, of order
	$$ G(j_n) \lesssim (\log{n})^{1+2\delta}\left(\frac{n}{\log{n}}\right)^{-\frac{2(r+1)}{2r+d/2}} = o(R_n(s)^2) = o(R_n(r)^2). $$
	(When $d\leq4$ this is trivial; for $d>4$, it necessitates the assumption on $s$.)
	
	Next, $n^{-2}S_{j_n}$ is of order
	$$ n^{-2}\sum_{l\leq j_n-1}2^{l(d-4)}(l\vee1)^{4\delta}. $$
	When $d\leq4$, this contributes at most a logarithmic factor in $n$ times $n^{-2}$, so this is clearly $o(R_n(s)^4)$ and $o(R_n(r)^4)$. When $d>4$, the final term dominates the sum and so the contribution is of order
	$$ (\log{n})^{4\delta - \frac{d-4}{2r + d/2}}n^{-\frac{4(r+1)}{2r+d/2}} = O(R_n(s)^4) = o(R_n(r)^4),$$
	again by the assumption on $s$.
	
	Finally, since $\mathrm{Var}(U_{j_n})\to0$ as $n\to\infty$, we know that
	$$ U_{j_n} = O_P\left(E_fU_{j_n}\right) = O_P\left(E_f\|K_j(f-\hat{f}_n)\|_{H^{-1,\delta}}^2)\right) = O_P\left(E_f\|f-\hat{f}_n\|_{H^{-1,\delta}}^2\right).$$
	As $\hat{f}_n$ converges at the rates $R_n(s)$ and $R_n(r)$ uniformly over $\mathcal{F}(s)$ and $\mathcal{F}(r)$ respectively, $U_{j_n}$ is of the correct order in probability in both cases. This concludes the proof.
\end{proof}

\begin{proof}[Proof of Theorem \ref{Thm: testing rate}]		
	For some sequence $L_n\to\infty$, to be defined below, and any $\omega\in \left\{-1;1\right\}^{\Z\cap\left[0,2^{L_n}\right)^d}$, we define for some $\epsilon>0$,
	\[  f_{n,\omega} \coloneqq 1 +  \epsilon2^{-L_n(r+d/2)}\sum_{k\in\Z\cap\left[0,2^{L_n}\right)^d}\omega_{k} \psi_{L_n,k}.\]
	Provided that $B>1$,
	\begin{align*}
		\norm{f_{n,\omega}}_{B^{r}_{2q}} &= 1 + 2^{L_nr} \left(\sum_{k\in\Z\cap\left[0,2^{L_n}\right)^d}|\langle f_{n,\omega}, \psi_{L_n,k}\rangle|^2\right)^{1/2}\\
		&= 1 + \epsilon 2^{L_nr} 2^{-L_n(r+d/2)} 2^{dL_n/2} \\
		&= 1 + \epsilon,
	\end{align*}
	ensuring that $ f_{n,\omega}$ is in the $\norm{\cdot}_{B^2_{2q}}$-Besov ball of radius $B$ for $\epsilon$ small enough. Also, $\int_{\mathbb{T}^d}f_{n,\omega}(t)dt=1$ and, as the tensor product wavelet basis is assumed to be $S-$regular (cf. Appendix \ref{Section: Wavelet Appendix}),
	\[ \norm{\sum_{k} |\psi_{L_n,k}|}_\infty \lesssim 2^{dL_n/2},\]
	for some constant depending on the basis only. Therefore,
	\[ \norm{f_{n,\omega}-1}_\infty \leq \epsilon c 2^{-rL_n},\]
	so that, for any  $M>1\geq m>0$, $f_{n,\omega}\in\mathcal{F}(r)$ for $n$ large enough (or $\epsilon$ small enough if $r=0$). Finally, for any $\rho_n=o\left(n^{-\frac{1+r}{2r+d/2}}\right)$, $f_{n,\omega}\in \tilde{\mathcal{F}}(r,\rho_n)$ if, for any $g\in \mathcal{F}(s)$, $W_2\left(f_{n,\omega}, g\right)\geq \rho_n$. By definition of $\mathcal{F}(r)$, $\mathcal{F}(s)$ and Proposition \ref{Prop: W-Besov comparison}, we have, for $n$ large enough
	\begin{align*}
		W_2\left(f_{n,\omega}, g\right)^2 &\gtrsim  \norm{f_{n,\omega} - g}^2_{B^{-1}_{2\infty}}\\
		&\geq 2^{-2L_n}\norm{\langle f_{n,\omega} - g, \psi_{L_n,\cdot}\rangle}_2^2\\
		&\geq  2^{-2L_n}\left[  \Bigg(\sum_{k=0}^{2^{L_nd}-1} |\langle f_{n,\omega} , \psi_{L_n,k}\rangle|^2\Bigg)^{1/2} - \Bigg(\sum_{k=0}^{2^{L_nd}-1}|\langle g, \psi_{L_n,k}\rangle|^2\Bigg)^{1/2}\right]^2\\
		&\geq  2^{-2L_n} \left[   \epsilon 2^{-L_nr} - B2^{-L_ns}\right]^2\\
		&\geq \frac{\epsilon^2}{2}2^{-2L_n(1+r)}.
	\end{align*}
	Therefore, if $L_n^*$ is such that $2^{-2L_n^*(1+r)}\asymp n^{-2\frac{1+r}{2r+d/2}}$, it is possible to find $L_n>L_n^*$ such that $\rho_n^2\leq \frac{\epsilon^2}{2}2^{-2L_n(1+r)} = o\left(n^{-2\frac{1+r}{2r+d/2}}\right)$. This choice ensures that, for any $\omega$, $f_{n,\omega}\in \tilde{\mathcal{F}}(r,\rho_n)$. Note also that the density $f_0\coloneqq 1$ naturally belongs to $\mathcal{F}(s)$.\\
	
	Re-index $\left\{-1;1\right\}^{\Z^d\cap\left[0,2^{L_n}\right)^d}$ as $\left\{\omega^{(i)}:\ i=1,\dots, 2^{2^{dL_n}} \right\}$ and denote by $P_i$ the distribution with Lebesgue density $f_i\coloneqq f_{n,\omega^{(i)}}$,  $Q:=2^{-2^{dL_n}} \sum_{i=1}^{2^{2^{L_nd}}}P_i$ and $P_0$ the distribution with density $f_0$. Then, with $\mu$ the Lebesgue measure and for any test $\Psi_n$,
	
	\begin{align*}
		\underset{f\in \Sigma_0}{\sup} \mathbb{E}_f \left[\Psi_n\right] +  \underset{f\in \Sigma(\rho_n)}{\sup} \mathbb{E}_f \left[1-\Psi_n\right] &\geq \mathbb{E}_{f_0} \left[\Psi_n\right] + 2^{-2^{dL_n}} \sum_{i=1}^{2^{2^{L_nd}}} \mathbb{E}_{f_i} \left[1-\Psi_n\right]\\
		&\geq \int \left(\Psi_n(x_1,\dots,x_n)+1-\Psi_n(x_1,\dots,x_n)\right)\\
		&\qquad\qquad\left(\prod_{j=1}^n f_0(x_j)\wedge  2^{-2^{dL_n}} \sum_{i=1}^{2^{2^{L_nd}}} \prod_{j=1}^n f_i(x_j)\right)d\mu^{\otimes n}(x_1,\dots,x_n)\\
		&= 1-\frac{1}{2}\norm{P_0^{\otimes n}- Q^{\otimes n} }_1\\
		&\geq 1-\frac{1}{2}\sqrt{\chi^2\left(Q^{\otimes n}, P_0^{\otimes n}\right)}. 
	\end{align*} where $\chi^2(Q,P)=\int (dP/dQ-1)^2 dQ$ if $P\ll Q$, $\chi^2(Q,P)=+\infty$ otherwise.
	Also, for any $1\leq \gamma,\kappa\leq 2^{2^{dL_n}}$, the orthonormality of the wavelet basis gives
	\begin{align*}
		&\int \frac{dP_\gamma^{\otimes n}}{dP_0^{\otimes n}}\frac{dP_\kappa^{\otimes n}}{dP_0^{\otimes n}}dP_0^{\otimes n} \\
		&= \prod_{i=1}^n \int_{\mathbb{T}^d} \left[1 + \epsilon2^{-L_n(r+d/2)}\sum_{k}\omega^{(\gamma)}_{k} \psi_{L_n,k}(x_i)\right]\left[1 + \epsilon2^{-L_n(r+d/2)}\sum_{k}\omega^{(\kappa)}_{k} \psi_{L_n,k}(x_i)\right]dx_i\\
		&= \left(1+\epsilon^2 2^{-L_n(2r+d)} \sum_{k}  \omega^{(\gamma)}_{k}  \omega^{(\kappa)}_{k} \right)^n.
	\end{align*}
	Then, for $\gamma_n=n\epsilon^22^{-L_n(2r+d)}\to0$ and $R_k, R_k'$ i.i.d. Rademacher random variables,
	\begin{align*}
		\chi^2\left(Q^{\otimes n}, P_0^{\otimes n}\right) &=  2^{-2^{dL_n}}\sum_{\gamma,\kappa} \left(1+\epsilon^2 2^{-L_n(2r+d)} \langle \omega^{(\gamma)}, \omega^{(\kappa)} \rangle\right)^n -1\\
		&\leq \mathbb{E}\left[\exp\Big(n\epsilon^22^{-L_n(2r+d)} \sum_{k} R_k R_k'\Big)\right]-1\\
		&= \mathbb{E}\left[\exp\Big(n\epsilon^22^{-L_n(2r+d)} \sum_{k} R_k\Big)\right]-1 \\
		&= \cosh(\gamma_n)^{2^{L_nd}} -1,
	\end{align*}
	where we used that $1+x\leq e^x$ for $x\in\R$ in the second line and that $R_kR_k'$ is distributed as $R_k$ in the third. Using that $\cosh(z)=1+z^2/2+\underset{|z|\to0}{o}(z^2)$ and $1+x\leq e^x$ once again, for any $\delta>0$,
	\[\left(\cosh(\gamma_n)\right)^{2^{dL_n}}-1=\left(1+\frac{\gamma_n^2}{2}(1+o(1))\right)^{2^{dL_n}}-1\leq \exp\left(\gamma_n^22^{dL_n-1}(1+o(1))\right)-1\leq \delta^2\]
	for $n$ large enough, since $\gamma_n^22^{dL_n}=o(1)$. We have proven that, for any $\beta<1$ and $\rho_n=o\left(\rho_n^*\right)$, 
	\[ \underset{n}{\lim\inf} \ \underset{\Psi_n}{\inf}\Big[\underset{f\in  \mathcal{F}(s)}{\sup} \mathbb{E}_f \left[\Psi_n\right] +  \underset{f\in \tilde{\mathcal{F}}(r,\rho_n)}{\sup} \mathbb{E}_f \left[1-\Psi_n\right] \Big] \geq \beta,\]
	which concludes the proof.
\end{proof}

	\section{Proofs for Section \ref{Section: Rd}}\label{Section: Rd proofs}

\begin{proof}[Proof of Proposition \ref{Prop: W1 loss decomposition}]
	As $f$ and $\tilde{f}_n$ have the same total mass, we may without loss of generality take the supremum over functions $h\in\mathrm{Lip}_1(\Rd)$ for which $h(0)=0$; observe that $x\mapsto\|x\|$ is an envelope for this function class. Since both $f$ and $\tilde{f}_n$ have finite first moments (almost surely), the wavelet expansion of any $h$ in this class converges in $L_1(f)$ and $L_1(\tilde{f}_n)$ and so 
	$$ \int_{\Rd}h(f-\tilde{f}_n) = \sum_{k\in\Zd}\langle h,\phi_k\rangle \langle f-\tilde{f}_n, \phi_k\rangle + \sum_{l\geq0}\sum_{k\in\Zd} \langle h,\psi_{lk} \rangle \langle f-\tilde{f}_n,\psi_{lk} \rangle. $$
	As the father wavelets $\phi_k$ are compactly supported in some interval about $k$,
	$$ |\langle h,\phi_k \rangle| \lesssim |h(k)| \leq \|k\| $$
	for some constant depending on the wavelet basis. Moreover, $h-K(h) = \sum_{l\geq0}\sum_{k\in\Zd}\langle h,\psi_{lk}\rangle\psi_{lk}$ is in a $B^1_{\infty\infty}$-ball of radius depending only on the wavelet basis, and so by (\ref{Eq: Besov wavelet coefficient bound}),
	$$ \sup_{k\in\Zd}|\langle h,\psi_{lk}\rangle| \lesssim 2^{-l\left(\frac{d}{2}+1\right)}. $$
	Plugging these uniform estimates for the wavelet coefficients of $h$ into the first equation gives the result.
\end{proof}

\begin{proof}[Proof of Theorem \ref{Thm: Rd estimation upper bound}]
	When $d=1$, the empirical measure achieves the stated rate (\cite{fournierRateConvergenceWasserstein2015}). Thus we assume $d\geq2$.
	
	The estimator we use is
	$$ \hat{f}_n := \sum_{\|k\|_{\infty}\leq \kappa_{-1n}} \hat{f}_{-1k}\phi_k + \sum_{l\leq l_n(s)}\sum_{\|k\|_{\infty}\leq\kappa_{ln}}\hat{f}_{lk}\psi_{lk}, $$
	where $\hat{f}_{lk}$ are empirical wavelet coefficients and the cutoffs $\kappa_{ln}, l_n(s)$ are chosen such that
	$$ 2^{l_n(s)} \simeq n^{\frac{1}{2s+d}}, \kappa_{-1n} = \kappa_{0n} \simeq (\log{n})^{\gamma}, \kappa_{ln} = 2^l\kappa_{0n},$$
	where $\gamma$ is to be chosen below.
	We then use the decomposition in Proposition \ref{Prop: W1 loss decomposition}, which we further split to obtain six terms:
	\begin{align*}
		W_1(f,\hat{f}_n) \lesssim& \sum_{\|k\|_{\infty}\leq \kappa_{-1n}} \|k\| |\hat{f}_{-1k}-f_{-1k}| + \sum_{\|k\|_{\infty}> \kappa_{-1n}} \|k\||f_{-1k}| \\
		&\ldots + \sum_{l< l_n(s)} 2^{-l\left(\frac{d}{2}+1\right)}\sum_{\|k\|_{\infty}\leq\kappa_{ln}}|\hat{f}_{lk}-f_{lk}| +  \sum_{l< l_n(s)} 2^{-l\left(\frac{d}{2}+1\right)}\sum_{\|k\|_{\infty}>\kappa_{ln}}|f_{lk}| \\
		&\ldots + \sum_{l\geq l_n(s)} 2^{-l\left(\frac{d}{2}+1\right)}\sum_{\|k\|_{\infty}\leq\kappa_{ln}}|f_{lk}| +  \sum_{l\geq l_n(s)} 2^{-l\left(\frac{d}{2}+1\right)}\sum_{\|k\|_{\infty}>\kappa_{ln}}|f_{lk}| \\
		&=: I + II + III + IV + V + VI.
	\end{align*}
	We first consider the bias terms $II,IV,VI$. For term $II$, we have that
	$$ \sum_{\|k\|_{\infty}> \kappa_{-1n}}\|k\||f_{-1k}| \leq \int_{\Rd}\sum_{\|k\|_{\infty}> \kappa_{-1n}}\|k\||\phi_k(x)|f(x)\,\ud x.$$
	Since each $\phi_k$ is compactly supported in some interval about $k$, and $\sum_{k\in\Zd}|\phi_k|$ is uniformly bounded on $\Rd$, we have that
	$$ \sum_{\|k\|_{\infty}> \kappa_{-1n}}\|k\||\phi_k(x)| \lesssim \|x\| $$
	for some constant depending on the wavelet basis. Moreover, the integrand is supported for all large enough $n$ in $([-\kappa_{-1n}/2,\kappa_{-1n}/2]^d)^c =: D_n$. Thus, for $n$ large enough,
	\begin{equation}\label{eq: bound bias mother wavelet 1} II \lesssim \int_{D_n}\|x\|f(x)\,\ud x \leq \mathcal{E}_{\alpha,\beta}(f) \kappa_{-1n}\exp\left(-\beta\left(\frac{\kappa_{-1n}}{2}\right)^{\alpha}\right). \end{equation}
	Since $\sum_{k\in\Zd}|\psi_{lk}| $ is uniformly bounded by a constant depending on the wavelet basis times $2^{ld/2}$, we analogously have
	\begin{equation}\label{eq: bound bias wavelet 2} \sum_{\|k\|_{\infty}>\kappa_{ln}} |f_{lk}| \lesssim 2^{ld/2}\mathcal{E}_{\alpha,\beta}(f)\exp\left(-\beta\left(\frac{\kappa_{0n}}{2}\right)^{\alpha}\right). \end{equation}
	Thus
	$$ IV + VI \lesssim \mathcal{E}_{\alpha,\beta}(f)\exp\left(-\beta\left(\frac{\kappa_{0n}}{2}\right)^{\alpha}\right). $$
	Choosing $\gamma>0$ sufficiently large depending on $\alpha,\beta$, these terms converge faster than $n^{-1/2}$.
	
	Next, we deal with the final bias term $V$. By Cauchy-Schwarz and the fact that $\|f\|_{B^s_{2q}}\leq B$,
	$$ \sum_{\|k\|_{\infty}\leq\kappa_{ln}} |f_{lk}| \leq \sqrt{\kappa_{ln}^d}\|f_{l\cdot}\|_2 \lesssim (\log{n})^{\gamma d/2}2^{l\left(\frac{d}{2} - s\right)},$$
	and so
	$$ V \lesssim \sum_{l\geq l_n(s)} 2^{-l(s+1)}(\log{n})^{\gamma d/2} \simeq (\log{n})^{\gamma d/2}2^{-l_n(s)(s+1)}$$
	which is of the correct order by the definition of $l_n(s)$.
	
	To bound the stochastic terms $I$ and $III$, we use the expectation bound Lemma \ref{Lemma: H-J inequality}, whose proof generalises naturally to the case of $\Rd$. We have for all $l\geq-1$ such that $2^{ld}\leq n$ and $k\in\Zd$ that
	$$ E_f|\hat{f}_{lk} - f_{lk}| \lesssim n^{-1/2},$$
	for some constant depending on $M$ and the wavelet basis. So
	$$ E_f(I) \lesssim (\kappa_{-1n})^{d+1}n^{-1/2} $$
	and
	$$ E_f(III) \lesssim (\log{n})^{\gamma d}n^{-1/2}\sum_{l< l_n(s)} 2^{l\left(\frac{d}{2}-1\right)}.$$
	When $d=2$, the sum contributes an extra $\log{n}$ factor as in the statement. For $d\geq3$, the final term of the sum dominates, and so
	$$ E_f(III) \lesssim (\log{n})^{\gamma d/2}n^{-\frac{s+1}{2s+d}} $$
	as stated.
\end{proof}

\begin{proof}[Proof of Theorem \ref{Thm: Rd adaptive estimator}]
Define the thresholds $\kappa_{-1n} = \kappa_{0n} \simeq (\log{n})^{\gamma}, \kappa_{ln} = 2^l\kappa_{0n}$ for $\gamma$ chosen as in Theorem \ref{Thm: Rd estimation upper bound}. As before, let $l_{\max}$ be such that $2^{l_{\max}} \simeq (n/\log{n})^{1/d}$; for $0\leq l\leq l_{\max}$, define the thresholds $\tau_l$ via
$$ \tau_l^2 = \tau^2\kappa_{ln}^d\frac{\log{n}}{n}, $$
where $\tau>0$ is to be chosen below. For any sequence $(a_k)_{k\in\Zd}$, set
$ \|a_{\cdot}\|_{2,\kappa_{ln}} := \left( \sum_{\|k\|_{\infty}\leq\kappa_{ln}} a_k^2 \right)^{1/2}. $
The thresholded estimator is then defined as
\begin{equation}\label{Eq: Rd thresholded estimator definition}
	\hat{f}_n = \sum_{\|k\|_{\infty}\leq \kappa_{-1n}}\hat{f}_{-1k}\phi_k + \sum_{l=0}^{l_{\max}}\ind_{\left\{ \|\hat{f}_{l\cdot}\|_{2,\kappa_{ln}}> \tau_l \right\}}\sum_{\|k\|_{\infty}\leq\kappa_{ln}} \hat{f}_{lk}\psi_{lk}.
\end{equation}

	We perform a decomposition of the risk similar to that in the previous proof:
	\begin{align*}
		W_1(f,\hat{f}_n) \lesssim& \sum_{\|k\|_{\infty}\leq \kappa_{-1n}} \|k\| |\hat{f}_{-1k}-f_{-1k}| + \sum_{\|k\|_{\infty}> \kappa_{-1n}} \|k\||f_{-1k}| \\
		&\ldots + \sum_{l\leq l_{\max}} 2^{-l\left(\frac{d}{2}+1\right)}\sum_{\|k\|_{\infty}\leq\kappa_{ln}}\left|\ind_{\left\{ \|\hat{f}_{l\cdot}\|_{2,\kappa_{ln}}> \tau_l \right\}}\hat{f}_{lk}-f_{lk}\right| +  \sum_{l\leq l_{\max}} 2^{-l\left(\frac{d}{2}+1\right)}\sum_{\|k\|_{\infty}>\kappa_{ln}}|f_{lk}| \\
		&\ldots + \sum_{l>l_{\max}} 2^{-l\left(\frac{d}{2}+1\right)}\sum_{\|k\|_{\infty}\leq\kappa_{ln}}|f_{lk}| +  \sum_{l> l_{\max}} 2^{-l\left(\frac{d}{2}+1\right)}\sum_{\|k\|_{\infty}>\kappa_{ln}}|f_{lk}| \\
		&=: I + II + III + IV + V + VI.
	\end{align*}
	We treat terms $I,II,IV$ and $VI$ identically to before. Term $V$ is also dealt with as in the previous proof, noting that for all $n$ sufficiently large, $2^{l_{\max}}>n^{1/(2s+d)}$. It remains to deal with term $III$; by Cauchy-Schwarz and the definition of $\kappa_{ln}$, we have that
	$$ III \lesssim (\log{n})^{\frac{\gamma d}{2}}\sum_{l=0}^{l_{\max}}2^{-l}\left\| \ind_{\left\{ \|\hat{f}_{l\cdot}\|_{2,\kappa_{ln}}> \tau_l \right\}}\hat{f}_{l\cdot} - f_{l\cdot} \right\|_{2,\kappa_{ln}},$$
	where the constant depends on $d$. By splitting this sum into two parts at $l_n(s)$, $2^{l_n(s)} \simeq B^{1/s}(n/\log{n})^{1/(2s+d)}$, one can bound it exactly as in the proof of Theorem \ref{Thm: Thresholded estimator}
\end{proof}

\begin{proof}[Proof of Theorem \ref{Thm: Rd adaptive confidence sets}]
	We first establish coverage. Define the thresholds $\kappa_{ln}$ as in the previous proof. Given $f\in\mathcal{G}(s)$, as in the proof of Theorems \ref{Thm: Thresholded estimator} and \ref{Thm: Rd adaptive estimator}, on an event of probability tending to 1, for all $l$ such that $l_n(r)\leq l\leq l_{\max}$, $\langle \hat{f}_n,\psi_{l\cdot}\rangle \equiv 0$. Note that $l_{\max}>j_n >l_n(s)>l_n(r)$. So on this event, by Cauchy-Schwarz,
	\begin{align*}
		\left(\sum_{l\geq j_n}2^{-l\left(\frac{d}{2}+1\right)}\sum_{\|k\|_{\infty}\leq\kappa_{ln}}|\langle f-\hat{f}_n,\psi_{lk}\rangle|\right)^2 &\lesssim (\log{n})^{\gamma d}\left(\sum_{l\geq j_n}2^{-l}\|\langle f,\psi_{l\cdot}\rangle\|_2\right)^2 \\
		&\lesssim (\log{n})^{\gamma d}B2^{-2j_n(r+1)} \\
		&\leq G_{j_n}
	\end{align*}
	for all $n$ sufficiently large, i.e. this quantity is $O_P(G_{j_n})$. The other term in $\tilde{W}^{(n)}(f,\hat{f}_n)^2$ is precisely $E^{(2)}_f V_{j_n}$; by Chebyshev's inequality we obtain condition (\ref{Eq: Rd honest coverage}).
	
	It remains to confirm the diameter conditions (\ref{Eq: Rd diameter shrinkage slow}) and (\ref{Eq: Rd diameter shrinkage fast}) with the rates $R_n(r),R_n(s)$ as given in the statement of the result. As the remainder term $\sqrt{r_n}$ converges up to a logarithmic factor at the rate $n^{-1/2}$, it is dominated by $\tilde{W}^{(n)}$ for diameter considerations. As observed previously, we may instead prove the diameter conditions for the $\tilde{W}^{(n)}$ distance with the augmented rates
	$$ \bar{R}_n(r) = (\log{n})^{\gamma d/2}\left(\frac{n}{\log{n}}\right)^{-\frac{r+1}{2r+d}}, \quad \bar{R}_n(s) = (\log{n})^{\gamma d/2}\left(\frac{n}{\log{n}}\right)^{-\frac{s+1}{2s+d}}. $$
	By the same argument as in the proof of Theorem \ref{Thm: d > 2 confidence sets}, the $\tilde{W}^{(n)}$-diameter of $C_n$ is bounded by a constant multiple of
	$$ (\log{n})^{\gamma d/4 + 1/2}n^{-1/2}\left(\sum_{l<j}2^{l(d-4)}\right)^{1/4} + \sqrt{V_{j_n}} + \sqrt{G_{j_n}} + n^{-1/2}. $$
	The final term is dominated by the first, and (using the condition on $s$ when $d>4$) $\sqrt{G_{j_n}} = O(\bar{R}_n(s)) = o(\bar{R}_n(r))$. One checks the first term is of the correct order as in Theorem \ref{Thm: d > 2 confidence sets}. Finally, since $\mathrm{Var}^{(2)}_f(V_{j_n})\to0$, we have that
	$$V_{j_n} = O_{P_f}\left(E_f V_{j_n}\right);$$
	as in the proof of Theorem \ref{Thm: Rd adaptive estimator}, this expectation is of order $\bar{R}_n(r)$ or $\bar{R}_n(s)$ when $f$ belongs to $\mathcal{G}(r)$ or $\mathcal{G}(s)$ respectively.
\end{proof}

\begin{proof}[Proof of Theorem \ref{Thm: testing rate Rd}]
	For some $\alpha'>\alpha$, $D>0$ and $\alpha(x)=\alpha'e^{-1/\left(\norm{x}_2-D\right)}\mathds{1}_{B(0,D)^c}(x)$, the density $f$ defined by \[f(x)\propto e^{-\beta\norm{x}_2^{\alpha(x)}}\]is such that $\mathbb{E}_f\Big[e^{\beta\norm{X}^\alpha}\Big]<+\infty$. Then, for $\sigma>0$, if $X\sim P_f$, $\sigma X$ has density $g: x\mapsto \sigma^{-d}f(\sigma^{-1} x)$ satisfying
	\[ \mathbb{E}_g\Big[e^{\beta\norm{X}^\alpha}\Big] = \mathbb{E}_f\Big[e^{\sigma^\alpha\beta\norm{X}^\alpha}\Big]<+\infty.\]
	Then, we verify that $f\in H_2^m(\R^d)\subset B_{2\infty}^{m}(\R^d)\subset B_{2q}^{s}(\R^d)$, for any $m\in N$ and $s<m$. Also, $\norm{g}_p=\sigma^{-d(1-1/p)}\norm{f}_p$ and, the moduli of continuity of $g$ satisifies, for $t>0$ and an integer $r>s$,
	\begin{align*}
		\omega_r(g,t,2) &\coloneqq \underset{0\leq \norm{h}\leq t}{\sup} \norm{\sum_{k=0}^r \binom{r}{k} (-1)^{r-k} \sigma^{-d} f(\sigma^{-1} \cdot + k\sigma^{-1} h)}_2 \\
		&=  \sigma^{-d} \underset{0\leq \norm{h}\leq \sigma^{-1} t}{\sup} \norm{\sum_{k=0}^r \binom{r}{k} (-1)^{r-k} \sigma^d f(\sigma^{-1} \cdot + k h)}_2\\
		&= \sigma^{-d/2}  \omega_r(f,\sigma^{-1} t,2).
	\end{align*}
	Therefore, with $\left|f\right|_{B^s_{pq}}\coloneqq \left[ \int_0^{\infty} \Big|\frac{\omega_r(f,t,p)}{t^s}\Big|^q\frac{dt}{t}\right]^{1/q}$, we have
	\begin{equation}\label{eqn: dilatation_besovnorm} \norm{g}_{B^{s}_{pq}} = \norm{g}_p + \left|g\right|_{B^s_{pq}} = \sigma^{-d(1-1/p)}\norm{f}_p +  \sigma^{-d(1-1/p)-s} \left|f\right|_{B^s_{pq}},\end{equation}
	so that $\norm{g}_{B_{2q}^s}\leq B$ for $\sigma$ large enough. Also, since $f\in L_\infty(\R^d)$, $g\leq M$ for $\sigma$ large enough. So, for some large $L$,  $g\in \mathcal{G}(s)$.

	For some sequence $L_n\to\infty$, and any $\omega\in \left\{-1;1\right\}^{\Z\cap\left[0,2^{L_n}\right)^d\times \mathcal{I}}$, we define for some $\epsilon>0$,
	\[ f^n_\omega = g + \epsilon2^{-L_n(r+d/2)}\sum_{k\in\Z\cap\left[0,2^{L_n}\right)^d, \iota\in\mathcal{I}} \omega_{k,\iota} \Psi^{\iota}_{L_nk}.\]
	Assuming that the scaling and mother wavelets functions are compactly supported (as assumed in Appendix \ref{Section: Wavelet Appendix}), the $\Psi^{\iota}_{L_nk}$, for $k\in\Z\cap\left[0,2^{L_n}\right)^d, \iota\in\mathcal{I}$, are supported on a compact set $K$ independent of $n$. Then, since
	\begin{align*}
		\norm{f^n_\omega}_{B^{r}_{2q}} &\leq \norm{g}_{B^{r}_{2q}} +  \epsilon2^{-L_n(r+d/2)} \norm{\sum_{k\in\Z\cap\left[0,2^{L_n}\right)^d, \iota\in\mathcal{I}} \omega_{k,\iota} \Psi^{\iota}_{L_nk}}_{B^{r}_{2q}}\\
		&\leq \norm{g}_{B^{r}_{2q}} +  C\epsilon,
	\end{align*}
	for some $C>0$ depending on $d$ only, reasoning as for \eqref{eqn: dilatation_besovnorm}, and since $B_{2q}^s\subset B_{2q}^r$, $f^n_\omega$ is the $\norm{\cdot}_{B_{2q}^r}$-Besov ball of radius $B$ for $\epsilon$ small enough and $\sigma$ large enough.
	Then, by assumption, $\int_{\R^d}f^n_\omega(t)dt=1$ and, since
	\[ \norm{2^{-L_n(r+d/2)}\sum_{k\in\R^d, \iota\in\mathcal{I}} \left|\Psi^{\iota}_{L_nk}\right|}_\infty \lesssim2^{-rL_n},\]
	$0< f^n_\omega \leq M$ for $n, \sigma$ large enough (or $\epsilon$ small enough if $r=0$). Indeed, $g$ is lower bounded by a some positive constant on $K$, So, $f^n_\omega$ actually is a density function.

	For these to belong to the alternative hypothesis, it remains to check that these are well separated from the null hypothesis. For any $h\in \mathcal{G}(s)$, the reversed triangular inequality gives
	
	\begin{align*}
		W_1\left(f^n_\omega, h\right)&\gtrsim  \norm{f^n_\omega - h}_{B^{-1}_{1\infty}}\\
		&\gtrsim  2^{-L_n(d/2+1)}   \sum_{k\in\Z\cap\left[0,2^{L_n}\right)^d, \iota\in\mathcal{I}}\left|\langle f^n_\omega - h, \Psi^{\iota}_{L_nk}\rangle\right|\\
		&\geq 2^{-L_n(d/2+1)}\left|  \sum_{k, \iota}|\langle f^n_\omega -g , \Psi^{\iota}_{L_nk}\rangle| - \sum_{k, \iota}|\langle h-g, \Psi^{\iota}_{L_nk}\rangle|\right|\\
		&= 2^{-L_n(d/2+1)} \left[  C 2^{-L_n(r-d/2)}- C'2^{-L_n(s-d/2)}\right]\\
		&\gtrsim2^{-L_n(1+r)},
	\end{align*} for constants independent of $n$. Above, we used that $s>r$ and that, for any $s>0$, $\mathcal{G}(s)\subset \left\{f: \norm{f}_{B^{s}_{1q}}\leq B'\right\}$ for some $B'>0$  according to Lemma \ref{Eq: exponential moment definition}..

%	Below, for any $h\in \mathcal{G}(s)$, $\left|\langle f^n_\omega - h, \Psi^{\iota}_{L_nk}\rangle\right|=\left| \epsilon \omega_{k,\iota} 2^{-L_n(r+d/2)} -\langle  h, \Psi^{\iota}_{L_nk}\rangle\right|\asymp 2^{-L_n(r+d/2)}$ since $\left|\langle  h, \Psi^{\iota}_{L_nk}\rangle\right|\lesssim 2^{-L_n s}$ and $s>r+d/2$. Cassel’s inequality \cite{MR89565} then implies
%	\begin{align*}
%		W_1\left(f^n_\omega, h\right)&\gtrsim  \norm{f^n_\omega - h}_{B^{-1}_{1\infty}}\\
%		&\gtrsim  2^{-L_n(d/2+1)}   \sum_{k\in\Z\cap\left[0,2^{L_n}\right)^d, \iota\in\mathcal{I}}\left|\langle f^n_\omega - h, \Psi^{\iota}_{L_nk}\rangle\right|\\
%		&\gtrsim  2^{-L_n}   \left(\sum_{k, \iota}\Big|\langle f^n_\omega - h, \Psi^{\iota}_{L_nk}\rangle\Big|^2\right)^{1/2}\\
%		&\gtrsim 2^{-L_n}\left[  \Bigg(\sum_{k, \iota}|\langle f^n_\omega -g , \Psi^{\iota}_{L_nk}\rangle|^2\Bigg)^{1/2} - \Bigg(\sum_{k, \iota}|\langle h-g, \Psi^{\iota}_{L_nk}\rangle|^2\Bigg)^{1/2}\right]\\
%		&\gtrsim 2^{-L_n} \left[  \epsilon \sqrt{2^d-1}2^{-L_nr}- 2B2^{-L_ns}\right]\\
%		&\gtrsim2^{-L_n(1+r)},
%	\end{align*} for some constant independent of $n$.
	The last inequality holds for $n$ large enough. Therefore, if $L_n^*$ is such that $2^{-L_n^*(1+r)}\asymp \xi_n$, it is possible to take $L_n>L_n^*$ such that $\rho_n\leq C'2^{-L_n(1+r)} = o(\xi_n)$, so that, for any $\omega$, $f^n_\omega \in \tilde{\mathcal{G}}(s,\rho_n)$.
	
	For $N_n=2^{2^{dL_n}\left(2^d-1\right)}$, let's index $\omega\in \left\{-1;1\right\}^{\Z\cap\left[0,2^{L_n}\right)^d\times \mathcal{I}}=\left\{w^{(m)}:\ m=1,\dots, N_n\right\}$ and denote $P_m=P_{f_{\omega^{(m)}}}$. Then, 
	\[ \underset{n}{\lim\inf} \ \underset{\Psi_n}{\inf}\Big[\underset{f\in H_0}{\sup} \mathbb{E}_f \left[\Psi_n\right] +  \underset{f\in H_1(r_n)}{\sup} \mathbb{E}_f \left[1-\Psi_n\right] \Big] \geq 1-\frac{1}{2}\sqrt{\chi^2\left(Q^{\otimes n}, P_0^{\otimes n}\right)},\]
	where $Q=N_n^{-1}\sum_{m=1}^{N_n} P_m$ and $P_0$ has density $g\in H_0$.
	Then, for any $1\leq m,q\leq N_n$, one has by properties of the wavelet basis, denoting $\nu_m=f_{\omega^{(m)}}-g$,
	\begin{align*}
		&\int \frac{dP_m^{\otimes n}}{dP_0^{\otimes n}}\frac{dP_q^{\otimes n}}{dP_0^{\otimes n}}dP_0^{\otimes n} \\
		&= \prod_{i=1}^n \int_{[0,1]^d} \left[g(x_i) + \epsilon2^{-L_n(r+d/2)}\sum_{k, \iota} \omega^{(m)}_{k,\iota} \Psi^{\iota}_{L_nk}(x_i)\right]\\
		&\qquad \left[g(x_i) + \epsilon2^{-L_n(r+d/2)}\sum_{k, \iota} \omega^{(q)}_{k,\iota} \Psi^{\iota}_{L_nk}(x_i)\right] g^{-1}(x_i)dx_i\\
		&= \left(1+ \int_{\R^d} \frac{\nu_m(x)\nu_q(x)}{g(x)} dx \right)^n.
	\end{align*}
	For $\sigma$ large enough, $g$ is constant on the compact support of $\nu_m$ and $\nu_q$, equal to $g(0)$. Hence, following the same arguments as above,
	\[
	\chi^2\left(Q^{\otimes n}, P_0^{\otimes n}\right) = (\cosh{\gamma_n})^{2^{dL_n}\left(2^d-1\right)}  -1,
	\]
	where $\gamma_n=n\epsilon^2g(0)^{-1}2^{-L_n(2r+d)}$, and for any $\delta>0$,
	$\chi^2\left(Q^{\otimes n}, P_0^{\otimes n}\right) \leq \delta^2$
	for $n$ large enough. This concludes the proof.
\end{proof}

\begin{lemma}\label{lemma: class inclusion}
Let $B\geq1,M>0,\alpha>0,\beta>0,L>0,1\leq q\leq \infty,$ and $s\geq0$. Then, there exists a constant $B'$, depending on the class parameters, the wavelet basis and the dimension $d$, such that
\[ \mathcal{G}_{s,2,q}(B,M;\alpha,\beta,L) \subset \mathcal{G}_{s,1,q}(B',M;\alpha,\beta,L). \]
\end{lemma}
\begin{proof}
Let $f\in \mathcal{G}_{s,2,q}(B,M;\alpha,\beta,L)$. All we have to prove is that \[\norm{f}_{B_{1q}^{s}}=\norm{\langle f, \phi_{\cdot}\rangle}_1+ \left(\sum_{l\geq 0} \big[2^{l(s-d/2)} \norm{\langle f, \psi_{l\cdot}\rangle}_1 \big]^q\right)^{1/q}\leq B',\] for some $B'$ as in the lemma. Let $\kappa>0$. Then, \[ \norm{\langle f, \phi_{\cdot}\rangle}_1 = \sum_{\norm{k}_\infty\leq \kappa} \left|\langle f,\phi_k\rangle\right|+\sum_{\norm{k}_\infty> \kappa} \left|\langle f,\phi_k\rangle\right|.\]
For the second term, the same arguments as the one used to obtain \eqref{eq: bound bias mother wavelet 1} give that it is bounded by $\mathcal{E}_{\alpha,\beta}(f) \exp\left(-\beta\left(\frac{\kappa}{2}\right)^{\alpha}\right)$, up to a constant depending on the wavelet basis. The first term is controlled via the Cauchy-Schwarz inequality
\[\sum_{\norm{k}_\infty\leq \kappa} \left|\langle f,\phi_k\rangle\right| \lesssim (2\kappa+1)^d \left(\sum_{\norm{k}_\infty\leq \kappa} \left|\langle f,\phi_k\rangle\right|^2\right)^{1/2}\leq (2\kappa+1)^d \norm{\langle f, \phi_{\cdot}\rangle}_2,\]
for a constant depending on $d$ only.

Next consider, for $l\geq0$, $\norm{\langle f, \psi_{l\cdot}\rangle}_1$. As before, letting $\kappa_l=2^{l/2}$, we have
\[ \norm{\langle f, \psi_{l\cdot}\rangle}_1 = \sum_{\norm{k}_\infty\leq \kappa_l} \left|\langle f,\psi_{lk}\rangle\right|+\sum_{\norm{k}_\infty> \kappa_l} \left|\langle f,\psi_{lk}\rangle\right|.\]
Arguing as with \eqref{eq: bound bias wavelet 2}, the second term is bounded by $2^{ld/2} \mathcal{E}_{\alpha,\beta}(f) \exp\left(-\beta\left(\frac{\kappa_l}{2}\right)^{\alpha}\right)$, up to a constant depending on the wavelet basis. The first term is controlled as above. Then, using the $l^q$ triangular inequality,
\begin{align*}
 \left(\sum_{l\geq 0} \big[2^{l(s-d/2)} \norm{\langle f, \psi_{l\cdot}\rangle}_1 \big]^q\right)^{1/q} &\lesssim  \left(\sum_{l\geq 0} 2^{ql(s-d/2)}\Big[ 2^{ld/2} \mathcal{E}_{\alpha,\beta}(f) \exp\left(-\beta\left(\frac{\kappa_l}{2}\right)^{\alpha}\right) +  (2\kappa_l+1)^d \norm{\langle f, \psi_{l\cdot}\rangle}_2 \Big]^q\right)^{1/q} \\
 &\lesssim \left(\sum_{l\geq 0}  \Big[2^{ls} \mathcal{E}_{\alpha,\beta}(f) \exp\left(-2^{-\alpha}\beta2^{l\alpha/2}\right) \Big]^q\right)^{1/q} + \left(\sum_{l\geq 0}\Big[ 2^{ls} \norm{\langle f, \psi_{l\cdot}\rangle}_2 \Big]^q\right)^{1/q},
\end{align*}
for a constants depending on the wavelet basis and $d$. The first term is upper bounded by \[\mathcal{E}_{\alpha,\beta}(f)  \left(\sum_{l\geq 0} 2^{qls}\exp\left(-q2^{-\alpha}\beta2^{l\alpha/2}\right) \right)^{1/q}\lesssim \mathcal{E}_{\alpha,\beta}(f),\] as the series converges.

In the end, following our assumptions on $\norm{f}_{B_{2q}^{s}}$, 
\begin{align*}
\norm{f}_{B_{1q}^{s}} &\lesssim (2\kappa+1)^d \norm{\langle f, \phi_{\cdot}\rangle}_2 + \left(\sum_{l\geq 0}\Big[ 2^{ls} \norm{\langle f, \psi_{l\cdot}\rangle}_2 \Big]^q\right)^{1/q} + \mathcal{E}_{\alpha,\beta}(f) \exp\left(-\beta\left(\frac{\kappa}{2}\right)^{\alpha}\right) + \mathcal{E}_{\alpha,\beta}(f)\\
&\lesssim B + \mathcal{E}_{\alpha,\beta}(f) \leq B+L,
\end{align*}
where the constants depend on the wavelet basis, $d$, the arbitrary $\kappa>0$ we took, $s$, $\alpha$, $\beta$ and $q$. 

\end{proof}

	\section{Wavelets and Besov Spaces}\label{Section: Wavelet Appendix}

Here we introduce the wavelet bases we use, and define the various norms and spaces used in our analysis.

\subsection{Wavelet Bases of $\Rd$ and $\Td$}

Let $S\in\mathbb{N}$. We begin with an $S$-regular wavelet basis of $L_2(\R)$ generated by scaling function $\Phi$ and wavelet function $\Psi$,
$$ \left\{ \Phi_k = \Phi(\cdot-k), \Psi_{lk} = 2^{l/2}\Psi(2^l(\cdot)-k): l\geq0, k\in\Z \right\}.$$
Concretely, we take sufficiently regular Daubechies wavelets: see \cite{meyerWaveletsOperators1993},\cite{daubechiesTenLecturesWavelets1992},\cite[Chapter 4]{gineMathematicalFoundationsInfiniteDimensional2015} for details. Such a wavelet basis has the following properties:
\begin{itemize}
	\item $\Phi,\Psi$ are in $C^S(\R)$, $\int_{\R}\Phi = 1$, and $\Psi$ is orthogonal to polynomials of degree $<S$.
	\item $\norm{\sum_k |\Phi_{k}|}_{\infty} \lesssim 1$, and$\norm{\sum_k |\Psi_{lk}|}_{\infty} \lesssim 2^{l/2}$ for a constant depending only on $\Psi$.
	\item Letting $V_j = \mathrm{span}(\Phi_k,\Psi_{lk}:l< j)$, for any $f\in V_j$ the following Bernstein estimate holds:
	$$ \norm{\nabla f}_p \lesssim 2^j\norm{f}_p,$$
	for a constant depending only on the wavelet basis. 
	\item $\Phi,\Psi$ are compactly supported.
%	\item For $\kappa(x,y)$ equal to either $\sum_{k\in\Z}\Phi_k(x)\Phi_k(y)$ or $\sum_k\Psi_{0k}(x)\Psi_{0k}(y)$, we have that
%	$$ \sup_{v\in\R}|\kappa(v,v-u)| \leq C\Omega(|u|) \quad \forall u\in\R, $$
%	where $C>0$ is a constant and $\Omega:[0,\infty)\to\R$ is bounded, integrable and such that $\int_{\R}|u|^S\Omega(|u|)\,\ud u <\infty$. \todo{Delete if we aren't discussing $L_1$ classes on $\Rd$...}
\end{itemize}

We then form a tensor product basis of $L_2(\Rd)$ as follows. Let $\mathcal{I}=\{0,1\}^d\setminus\{0\}$. Define
$$ \phi(x) = \Phi(x_1)\cdots\Phi(x_d), \quad x\in\Rd$$
and, writing $\Psi^0=\Phi,\Psi_1=\Psi$,
$$ \psi^{\iota} = \Psi^{\iota_1}(x_1)\cdots\Psi^{\iota_d}(x_d), \quad \iota\in\mathcal{I}.$$
Then (\cite[Section 4.3.6]{gineMathematicalFoundationsInfiniteDimensional2015})
$$ \left\{ \phi_k = \phi(\cdot-k), \psi^{\iota}_{lk} = 2^{ld/2}\psi^{\iota}(2^lx-k) : l\geq0,k\in\Zd,\iota\in\mathcal{I} \right\} $$
defines a wavelet basis of $L_2(\Rd)$. We omit $\iota$ from our notation and simply write $\psi_{lk}$ with $k$ now implicitly taking values in $\Zd\times\mathcal{I}$; any sum over $k$ is to be understood as over all $\iota\in\mathcal{I}$ as well. 
\begin{enumerate}[1)]
	\item $\phi,\psi$ are in $C^S(\Td)$, $\int_{\Rd}\phi = 1$, and $\psi$ is orthogonal to polynomials of degree $<S$.
	\item $\norm{\sum_k |\phi_{k}|}_{\infty} \lesssim 1$, and$\norm{\sum_k |\psi_{lk}|}_{\infty} \lesssim 2^{ld/2}$ for a constant depending only on $\psi$. 
	\item $\phi,\psi$ are compactly supported.
%	\item For $\kappa(x,y)$ equal to either $\sum_{k\in\Zd}\phi_k(x)\phi_k(y)$ or $\sum_k\psi_{0k}(x)\psi_{0k}(y)$, we have that
%	$$ \sup_{v\in\R^d}|\kappa(v,v-u)| \leq C\Omega(\|u\|) \quad \forall u\in\Rd, $$
%	where $C>0$ is a constant and $\Omega:[0,\infty)\to\R$ is bounded, integrable and such that $\int_{\Rd}\|u\|^S\Omega(\|u\|)\,\ud u <\infty$.  \todo{Same again}
\end{enumerate}
These properties follow elementarily from the previously stated properties of $\Phi$ and $\Psi$. Property 3) is used crucially in our analysis on $\Rd$. Notably, this precludes certain common choices of wavelet basis, such as the Meyer basis.
% \todo{Some comment about property 4)? (Assuming we keep it...)}

These properties imply the following relationship between $L_p$-norms of functions and the $\ell_p$-norms of their wavelet coefficients (by an abuse of notation we denote both of these norms by $\|\cdot\|_p$).
\begin{lemma}\label{Lemma: Lp-lp norm wavelet equivalence}
	For any $l\geq0$, any $p\in[1,\infty]$ and any $c\in\R^{\Zd}$, we have that
	$$ \norm{\sum_{k\in\Zd} c_k\psi_{lk}}_p \simeq 2^{ld(1/2 - 1/p)}\|c\|_p, $$
	where the constants depend on $\psi$ and $p$ only.
\end{lemma}

When working on $\Td$, we use the tensor product wavelet basis induced by the periodisations of $\Phi,\Psi$; see \cite[Section 4.3.4]{gineMathematicalFoundationsInfiniteDimensional2015} for details. This produces a basis of $L_2(\Td)$ with the following properties:
\begin{enumerate}[1)]
	\item $\psi(x) = \prod_{i=1}^d\psi^{(i)}(x_i)$ for some univariate functions $\psi^{(i)}$.
	\item Setting $\psi_{lk}(\cdot) = 2^{ld/2}\psi(\cdot - 2^{-l}k)$ for $l\geq0,k\in\Zd\cap[0,2^l)^d$, the set
	$$ \left\{ \phi,\psi_{lk} : l\geq0,k\in\Zd\cap[0,2^l)^d \right\} $$
	forms an orthonormal basis of $L_2(\Td)$. By an abuse of notation, we re-index in $k$ such that $k\in\Z$ varies between $0\leq k<2^{ld}$.
	\item $\psi$ is in $C^S(\Td)$, and is orthogonal to polynomials of degree $<S$.
	\item $\norm{\sum_k |\psi_{lk}|}_{\infty} \lesssim 2^{ld/2}$, for a constant depending only on $\psi$.
	\item Letting $V_j = \mathrm{span}(\phi,\psi_{lk}:l< j)$, for any $f\in V_j$ the following Bernstein estimate holds:
	$$ \norm{\nabla f}_p \lesssim 2^j\norm{f}_p,$$
	for a constant depending only on the wavelet basis. 
\end{enumerate}
Again, these are basic consequences of properties of $\Phi,\Psi$, and enable the proof of Proposition \ref{Prop: W-Besov comparison}; compare to Appendix C of \cite{weedEstimationSmoothDensities2019}.

\subsection{Besov Spaces}

In this section, we let $(\phi_k,\psi_{lk})$ denote either the $S$-regular tensor product Daubechies wavelet basis of $L_2(\Rd)$, or the $S$-regular tensor product periodised Daubechies wavelet basis of $L_2(\Td)$. It should be understood that any summation is over the full range of indices, for example $\sum_k\psi_{lk}$ denotes $\sum_{k\in\Zd}\psi_{lk}$ in the $\Rd$ case and $\sum_{k=0}^{2^{ld}-1}\psi_{lk}$ in the $\Td$ case.
We further let $\mathcal{D}$ be either the class of tempered distributions on $\Rd$, or the class of periodic tempered distributions on $\Td$. 

Let $1\leq p\leq \infty$, $1\leq q\leq \infty$, $s\in\mathbb{R}, s<S$. For $f\in\mathcal{D}$, we define the Besov norm
\begin{equation}\label{Eq: Besov norm definition}
	\|f\|_{B^s_{pq}} = \norm{\langle f,\phi_{\cdot}\rangle}_p + \left(\sum_{l\geq0}\left[2^{ls}2^{ld\left(\frac{1}{2}-\frac{1}{p}\right)}\norm{\langle f,\psi_{l\cdot}\rangle}_p\right]^q\right)^{1/q},
\end{equation}
where $\|\cdot\|_{p}$ is the $\ell_p$-norm. When $q=\infty$, the norm is defined as
\begin{equation}\label{Eq: Besov norm definition, q=infty}
	\|f\|_{B^s_{p\infty}} = \norm{\langle f,\phi_{\cdot}\rangle}_p + \sup_{l\geq0}2^{ls}2^{ld\left(\frac{1}{2}-\frac{1}{p}\right)}\norm{\langle f,\psi_{l\cdot}\rangle}_p.
\end{equation}
We then define the corresponding Besov space $B^s_{pq}$ as
\begin{equation}\label{Eq: Besov space definition}
	B^s_{pq} = \left\{f\in\mathcal{D}: \|f\|_{B^s_{pq}}<\infty \right\}.
\end{equation}
We will write $B^s_{pq}(\Rd)$ or $B^s_{pq}(\Td)$ to remove any ambiguity over the choice of domain, whenever it arises.

The definition of $B^s_{pq}$ is independent of the wavelet basis used, that is, using a different (sufficiently regular) basis in the definition (\ref{Eq: Besov norm definition}) produces an equivalent norm. Moreover, using a $C^{\infty}$ basis such as the Meyer basis enables us to define $B^s_{pq}$ concurrently for all $s\in\R$.

\subsection{The Case of the Unit Cube}\label{Subsection: unit cube}

We can also define a `boundary-corrected' wavelet basis of $L_2([0,1]^d)$ based on $\Phi,\Psi$ as in \cite{cohenWaveletsIntervalFast1993}; see also \cite[Section 4.3.5]{gineMathematicalFoundationsInfiniteDimensional2015}. Such a basis possesses completely analogous properties to properties 1)-5) of the periodised basis of $L_2(\Td)$; moreover, all Besov spaces defined on $\Td$ are defined on the unit cube $[0,1]^d$ by replacing the periodised wavelet basis with the boundary-corrected wavelet basis (as used in \cite{weedSharpAsymptoticFinitesample2019}). Thus all of our results for $\Td$ hold also for the case of $[0,1]^d$.

%	\listoftodos
	
	\bibliographystyle{apalike}
	\bibliography{wassersteinminimaxtheory5}
	
\end{document}